\documentclass[12pt,a4paper]{article}

\usepackage{amsmath}
\usepackage{amssymb}
\usepackage{mathrsfs}
\usepackage[colorlinks]{hyperref}
\usepackage{here}
\usepackage{xcolor,graphicx}
\graphicspath{{./input/}}

\newtheorem{theorem}{Theorem}[section]
\newtheorem{remark}[theorem]{Remark}
\newcommand{\de}{{\rm d}}
\newcommand{\bX}{{\mathbf{X}}}
\newcommand{\bx}{{\boldsymbol{x}}}
\newcommand{\bzeta}{{\boldsymbol{\zeta}}}
\newcommand{\bz}{{\mathbf{z}}}

\newcommand{\beq}{\begin{equation}}
\newcommand{\eeq}{\end{equation}}

\renewcommand{\contentsname}{}

\begin{document}

\title{
	\vspace{-2cm}\emph{Koopmon} trajectories\\ in nonadiabatic quantum-classical dynamics
}

\author{
	Werner Bauer$^1$,
	Paul Bergold$^1$\footnote{Corresponding author ({\tt p.bergold@surrey.ac.uk})}\,,
	Fran\c{c}ois Gay-Balmaz$^2$,
	Cesare Tronci$^{1}$\smallskip\\
	\small
	\it $^1$Department of Mathematics, University of Surrey, Guildford, United Kingdom\\
	\small
	\it $^2$Division of Mathematical Sciences, Nanyang Technological University, Singapore
}
\date{\vspace{-1.8cm}}

\maketitle
\begin{abstract}\footnotesize
	In order to alleviate the computational costs of fully quantum nonadiabatic dynamics, we present a mixed quantum-classical (MQC) particle method based on the theory of Koopman wavefunctions.
	Although conventional MQC models often suffer from consistency issues such as the violation of Heisenberg's principle, we overcame these difficulties by blending Koopman's classical mechanics on Hilbert spaces with methods in symplectic geometry.
	The resulting continuum model enjoys both a variational and a Hamiltonian structure, while its nonlinear character calls for suitable closures.
	Benefiting from the underlying action principle, here we apply a regularization technique previously developed within our team.
	This step allows for a singular solution ansatz which introduces the trajectories of computational particles -- the \emph{koopmons} -- sampling the Lagrangian classical paths in phase space.
	In the case of Tully's nonadiabatic problems, the method reproduces the results of fully quantum simulations with levels of accuracy that are not achieved by standard MQC Ehrenfest simulations.
	In addition, the koopmon method is computationally advantageous over similar fully quantum approaches, which are also considered in our study.
	As a further step, we probe the limits of the method by considering the Rabi problem in both the ultrastrong and the deep strong coupling regimes, where MQC treatments appear hardly applicable.
	In this case, the method succeeds in reproducing parts of the fully quantum results.
\end{abstract}
\vspace{-12mm}

{\footnotesize
\contentsname
\tableofcontents
}
\section{Introduction}
\subsection{Trajectory-based schemes in molecular dynamics}
In quantum molecular dynamics, several trajectory-based algorithms have been proposed over the years to overcome the computational challenges of fully quantum simulations.
While many of these algorithms are traditionally based on the remarkable properties of Gaussian wavepackets \cite{Heller,Lasser}, a considerable effort has also been devoted to exploiting Madelung's hydrodynamics \cite{Agostini,Garashchuk10,Kendrick,Wyatt}, which has the great advantage of restoring the concept of Lagrangian trajectory in quantum mechanics.
Severe difficulties, however, challenge both approaches.
For example, in the first case, the appearance of generally singular matrices often leads to uncontrolled approximations with the purpose of eliminating the degeneracy \cite{ShCh04}.
As a result \cite{JoIz18}, ``the existent formalisms using classical frozen-width Gaussian motion do not conserve the total energy'' and also struggle with momentum balance and conservation of probability \cite{HaWeTrBNMa01}.
Indeed, these limitations generally affect multiconfiguration methods and other approaches based on basis-set expansions; see \cite{JoIz18} for further details.

Madelung-based particle methods, on the other hand, often present striking difficulties related to the so-called \emph{quantum potential}, which renders the wave-like nature of quantum systems and leads to the emergence of high-frequency interference patterns posing severe computational challenges \cite{ZhMa03}.
Thus, in a typical molecular problem involving the correlation dynamics of interacting nuclei and electrons, one is tempted to simply discard the quantum potential contributions associated to the nuclei while retaining the quantum character of the electrons.

Arising naturally within the Born--Oppenheimer adiabatic approximation, this mixed quantum~-classical (MQC) description represents an attractive avenue also in the more general context of nonadiabatic molecular dynamics.
Yet, new problems emerge in extending the MQC method.
For example, even if one discards the nuclear quantum potential, the need to retain correlation effects leads to dealing with the so-called \emph{quantum geometric tensor} \cite{FoHoTr19, Provost80}, whose computational complexity is comparable to that of the quantum potential itself.
For these reasons, MQC models continue to represent a subject of research in computational chemistry.
Other areas in which the MQC description is attractive range from the theory of gravity \cite{BoDi21} to spintronics \cite{Petrovic}, and solid-state physics \cite{BACaCaVe09,HuHeMa17}.
A possible realistic realization of an MQC system is represented by the molecular reaction dynamics occurring within a liquid solvent, for which the classical fluid description applies with great accuracy \cite{Bo11, GBTr-fluid}.
 
\subsection{Consistency of MQC theories and the Ehrenfest model\label{sec:Ehrenfest}}
The fact that the formulation of MQC models represents a formidable challenge should not surprise us.
The correlation dynamics of interacting quantum and classical systems has represented a fundamental open question ever since the early days in the history of quantum mechanics, and goes back to the measurement problem.
Several remarkable attempts were made over the decades.
Within the chemistry community, a popular model is known under the name of \emph{quantum-classical Liouville equation} \cite{Aleksandrov,Schutte,Gerasimenko,Kapral}
\beq\label{QCLE}
	i\hbar\frac{\partial\widehat{\cal P}}{\partial t}+i\hbar\operatorname{div}\!\big(\widehat{\cal P}\bX_{\widehat{H}}\big)^{\sf H}
	=\big[\widehat{H},\widehat{\cal P}\big].
\eeq
Here, $\widehat{H}(q,p)$ is a Hamiltonian (Hermitian) operator on the quantum Hilbert space $\mathscr{H}$ parameterized by the classical phase-space coordinates and $\hbar$ denotes the reduced Planck constant.
Also, $\widehat{\cal P}$ is a sufficiently smooth distribution taking values in the space $\operatorname{Her}(\mathscr{H})$ of Hermitian operators on $\mathscr{H}$.
The classical Liouville density and the quantum density matrix are given as
\begin{align*}
	f
	=\operatorname{Tr}\widehat{\cal P}
	\qquad\text{and}\qquad
	\hat\rho
	=\int\!\widehat{\cal P}\,\de q\de p,	
\end{align*}
respectively.
In addition, we have introduced the operator-valued Hamiltonian vector field $\bX_{\widehat{H}}=(\partial_p{\widehat{H}},-\partial_q{\widehat{H}})$, while the superscript $\sf H$ denotes the projection on the Hermitian part so that $\widehat{A}^{\,\sf H}=(\widehat{A}+\widehat{A}^{\,\dagger})/2$.
Obtained as a limit case of fully quantum dynamics, equation \eqref{QCLE} generally fails to retain positivity of the classical density and the quantum density matrix, thereby invalidating the uncertainty principle \cite{AgCi07}.
This violation may or may not appear explicitly depending on the particular problem under consideration.
Analogously, the \emph{surface hopping} scheme in computational chemistry \cite{Tully} suffers from similar problems \cite{Bondarenko}, while other consistency issues also affect alternative MQC descriptions.

A list of consistency criteria to be satisfied by MQC models was presented in \cite{Traschen} and here we present five of those as follows:
1) the classical system is identified by a probability density at all times;
2) the quantum system is identified by a positive-semidefinite density operator $\hat\rho$ at all times;
3) the model is equivariant under both quantum unitary transformations and classical canonical transformations;
4) in the absence of a quantum-classical interaction potential, the model reduces to uncoupled quantum and classical dynamics;
5) in the presence of an interaction potential, the {\it quantum purity} $\|\hat\rho\|^2$ has nontrivial dynamics (decoherence property).
While the last item is absent in the list from \cite{Traschen}, it is included here to eliminate common mean-field descriptions which fail in reproducing decoherence.
We note that the list from \cite{Traschen} also included an item requiring linearity of the model equations, although this requirement is eventually questioned in the same reference.
Indeed, while there are reasons for the linear evolution of quantum wavefunctions \cite{Gisin,Polchinski}, it is not yet clear whether similar linearity criteria should also apply to MQC models.
If one insists on linearity as a requirement in MQC dynamics, one is naturally led to consider Lindblad-type structures for the irreversible evolution of MQC systems \cite{Diosi}.

If, on the other hand, we allow for nonlinear evolution equations, the only established reversible model satisfying all the above criteria and enjoying a Hamiltonian structure is the \emph{Ehrenfest model}
\beq\label{Ehrenfest}
	i\hbar\frac{\partial\widehat{\cal P}}{\partial t}
	+i\hbar\operatorname{div}\!\big(\widehat{\cal P}\langle\bX_{\widehat{H}}\rangle\big)
	=\big[\widehat{H},\widehat{\cal P}\big].
\eeq
Here, we have used the notation $\langle\widehat{A}\rangle=\operatorname{Tr}(\widehat{\cal P}\widehat{A})/\operatorname{Tr}\widehat{\cal P}$.
We observe that $\widehat{\cal P}$ undergoes a unitary evolution (due to the commutator) while being Lie-transported (due to the divergence) by the phase-space paths generated by the vector field $\langle\bX_{\widehat{H}}\rangle$.
Consequently, $\widehat{\cal P}$ is also preserved to be Hermitian and positive semidefinite.
In computational chemistry, the model \eqref{Ehrenfest} is used extensively in different variants \cite{Vanicek} due to the computational simplicity of its trajectory scheme.
This is obtained by first rewriting \eqref{Ehrenfest} as ${\partial_t}f+\operatorname{div}(f\langle\bX_{\widehat{H}}\rangle)=0$ and $i\hbar\partial_t\,\widehat{\!\mathscr{P}}+i\hbar\langle\bX_{\widehat{H}}\rangle\cdot\nabla\,\widehat{\!\mathscr{P}}=[\widehat{H},\,\widehat{\!\mathscr{P}}]$, where $\,\widehat{\!\mathscr{P}}=\widehat{\cal P}/f$ and $\langle\bX_{\widehat{H}}\rangle=\operatorname{Tr}(\,\widehat{\!\mathscr{P}}\bX_{\widehat{H}})$.
Then, the first equation for $f$ is solved by the point-particle solution ansatz
\beq\label{singansatz1}
	f(\bz,t)
	=\sum_{a=1}^Nw_a\delta(\bz-\bzeta_a(t)),
\eeq
with $\bz=(q,p),\,\bzeta_a(t)=(q_a(t),p_a(t)),\,\dot{\bzeta}_a=\langle\bX_{\widehat{H}}\rangle|_{\bz=\bzeta_a}$, and positive weights $w_a>0$ such that $\sum_aw_a=1$.
Here, the quantity $\langle\bX_{\widehat{H}}\rangle|_{\bz=\bzeta_a}$ requires evaluating ${\hat\varrho_a(t):=\,\widehat{\!\mathscr{P}}(\bzeta_a(t),t)}$ at all times, and this can be done by integrating the equation \eqref{Ehrenfest} for $\widehat{\cal P}=f\,\widehat{\!\mathscr{P}}$, so that $i\hbar\dot{\hat\varrho}_a=[\widehat{H}_a,\hat\varrho_a]$.
The resulting \emph{multi-trajectory Ehrenfest system}
\beq\label{MFeqs}
	\dot{q}_a
	=\partial_{p_a\!}\langle\hat\varrho_a|\widehat{H}_a\rangle,\qquad\
	\dot{p}_a
	=-\partial_{q_a\!}\langle\hat\varrho_a|\widehat{H}_a\rangle,\qquad\
	i\hbar\dot{\hat\varrho}_a
	=[\widehat{H}_a,\hat\varrho_a],
\eeq
with $\widehat{H}_a=\widehat{H}(\bzeta_a)$ and $\langle\hat\varrho_a|\widehat{H}_a\rangle=\operatorname{Tr}(\hat\varrho_a\widehat{H}_a)$, inherits all the consistency properties listed above \cite{Alonso}.
In particular, upon using $\hat\rho=\int\!f\,\widehat{\!\mathscr{P}}\de q\de p$, the quantum purity is expressed as $\operatorname{Tr}(\hat\rho^2)=\|\hat\rho\|^2=\|\sum_a w_a\hat\varrho_a\|^2$.
This scheme is also cost-effective from a computational perspective \cite{FaJiSp18}, although it is known to fail in reproducing accurate levels of decoherence.
For this and related reasons, much of the community has been looking for MQC approaches beyond the Ehrenfest model.
At present, MQC modeling comprises an entire field of research in computational chemistry, as demonstrated by popular recent reviews \cite{CrBa18}.
While most of the available MQC algorithms are based on ad-hoc methods and uncontrolled approximations, here we want to propose an alternative approach that is based on modern mathematical methods combining the geometry of Hamiltonian and Lagrangian systems with the Lie-group actions involved in their time evolution.

\subsection{Beyond the Ehrenfest model}
Much inspired by the work of George Sudarshan \cite{Marmo,Sudarshan}, we recently introduced a new nonlinear MQC theory that satisfies all the five consistency criteria listed previously and, at the same time, enjoys a Hamiltonian and variational structure.
The formulation of this model involves two main steps that exploit and blend unconventional methods from symplectic geometry with the theory of Koopman wavefunctions in classical mechanics \cite{Koopman}.
The latter are simply square-integrable functions $\chi\in L^2(T^*Q)$
on phase space such that the classical Liouville density is expressed as $f=|\chi|^2$.
Here, $T^*Q$ denotes the classical phase space, which is identified with $\mathbb{R}^2$ in this paper for the sake of simplicity. 

The \emph{Koopman-von Neumann equation} (KvN) of classical mechanics is written in terms of the canonical Poisson bracket as $i\hbar\partial_t\chi=\{i\hbar H,\chi\}$. 
The latter means that the classical evolution of Koopman wavefunctions is given as $\chi(t)=\chi_0\circ{\boldsymbol\eta}_t^{-1}$, that is the composition of the initial condition $\chi_0$ by the inverse of the Lagrangian map ${\boldsymbol\eta}_t$ generated by the vector field $\bX_H$.
We remark that the composition map is also known as the \emph{Koopman operator} and has recently attracted a great deal of attention in dynamical systems and data science \cite{Klus16,Li17,Mezic22}.

Importantly, the KvN equation identifies a linear unitary flow analogous to Schr\"odinger's quantum mechanics.
Due to the relation $f=|\chi|^2$, the KvN equation may be modified to include an arbitrary local phase factor, that is, $i\hbar\partial_t\chi=\{i\hbar H,\chi\}+\varphi\chi$ with $\varphi(q,p)$ being an arbitrary function.
Since both quantum and classical dynamics can be formulated in terms of wavefunctions, the idea of considering their tensor-product space is rather intuitive and identifies the starting point of Sudarshan's proposal.
However, writing the dynamics for a MQC wavefunction $\Upsilon(q,p,x)$ on such a tensor-product space (here, $x$ is the quantum coordinate) is far from intuitive and we provide a summary of our derivation in the following section.

Similarly to the Ehrenfest model \eqref{Ehrenfest}, the model proposed in \cite{GBTr22,GBTr21} formulates the dynamics of an operator-valued distribution $\widehat{\cal P}$ whose matrix-element representation is expressed in terms of the aforementioned MQC wavefunction as $\widehat{\cal P}(q,p,x,x')=\Upsilon(q,p,x){\Upsilon}(q,p,x')^*$, or, in finite dimensions, $\widehat{\cal P}(q,p)=\Upsilon(q,p)\Upsilon(q,p)^\dagger$.
Upon recalling $f=\operatorname{Tr}\widehat{\cal P}$ and introducing the quantity $\boldsymbol{\widehat{\Sigma}}=if^{-1}[\widehat{\cal P},\bX_{\widehat{\cal P}}]$, the MQC model proposed in \cite{GBTr22,GBTr21} reads
\beq\label{HybEq1}
	i\hbar\frac{\partial\widehat{\cal P}}{\partial t}+i\hbar\operatorname{div}\!\big(\widehat{\cal P}\boldsymbol{\cal X}\big)
	=\big[\widehat{\cal H},\widehat{\cal P}\big],
\eeq
with
\beq\label{HybEq2}
	\boldsymbol{\cal X}
	=\langle\bX_{\widehat{H}}\rangle+\frac\hbar{2f}\operatorname{Tr}\!\big(\bX_{{\widehat{H}}}\cdot\nabla\boldsymbol{\widehat{\Sigma}}-\boldsymbol{\widehat{\Sigma}}\cdot\nabla\bX_{{\widehat{H}}}\big),
\eeq
and
\beq\label{HybEq3}
	\widehat{\mathcal{H}}
	=\widehat{H}+\frac{i\hbar}f\Big(\{\widehat{\cal P},{\widehat{H}}\}+\{\widehat{H},\widehat{\cal P}\}+\frac{1}{2}[\{\ln f,\widehat{H}\},\widehat{\cal P}]\Big).
\eeq
We observe that both the vector field \eqref{HybEq2} and the effective Hamiltonian \eqref{HybEq3} are $\hbar$-corrections of the corresponding quantities appearing in the Ehrenfest model \eqref{Ehrenfest}, that is $\langle\bX_{\widehat{H}}\rangle$ and $\widehat{H}$, respectively.
Recovering Born--Oppenheimer molecular dynamics in the adiabatic approximation \cite{BeTr23}, this model was recently extended to consider MQC spin systems \cite{GBTr23} and solute-solvent interaction in fluid solvents \cite{GBTr-fluid}.
The corresponding Heisenberg picture is also available in \cite{DMCTr24}.

As discussed in \cite{TrGB23}, despite the formidable look of its equations, this system does not involve differentials of order higher than two.
Nevertheless, the high level of complexity of \eqref{HybEq1}-\eqref{HybEq3} makes these equations hardly tractable by standard methods.
In order to formulate and implement a numerical approach, this paper presents a particle method based on the variational structure underlying the continuum model.
In particular, a suitable regularization is applied to the action principle in such a way to allow for singular solutions that are otherwise unavailable.
The resulting point trajectories are used to sample the Lagrangian paths generated by the vector field $\boldsymbol{\cal X}$ in \eqref{HybEq2}.
As a result of this variational approach, the method inherits a Hamiltonian structure, which ensures the energy and momentum balance for all parameter regimes.
In addition, quantum decoherence is retained with accuracy levels that seem unachievable by standard Ehrenfest dynamics.
As expected in MQC models \cite{Agostini}, the resulting trajectories are all coupled together, thereby raising the level of complexity beyond Ehrenfest dynamics.
Yet, when compared to analogous trajectory-based fully quantum simulations, we observe that the proposed MQC model possesses substantial computational advantages.
First formulated in \cite{FoHoTr19,TrGB23}, this geometric mechanics approach reveals the Lagrangian paths of MQC dynamics.
Given the prominent role of Koopman wavefunctions in the formulation of the continuum model \eqref{HybEq1}-\eqref{HybEq3} underlying the numerical method, we will call its computational particles \emph{koopmons}.
Analogously, the method itself will be referred to as the \emph{koopmon method}.

\paragraph{Plan of the paper.}
The paper is organized as follows:
In Section~\ref{sec:The mixed quantum-classical koopmon method}, we introduce the koopmon method.
Starting from a brief review of the derivation and variational structure of the continuum model \eqref{HybEq1}-\eqref{HybEq3}, we demonstrate how the particle method emerges form of regularization at the level of the underlying variational principle.
Afterwards, in Section~\ref{sec:Numerical implementation}, we discuss the implementation of the koopmon method, which was used to produce the numerical experiments.
In particular, we discuss its connection to the quantum \emph{bohmion method}, and provide important details regarding the quantum solver and the output visualization.
In Section~\ref{sec:Results for the Tully models}, we present the results for the Tully models (I-III), which are well-known benchmark models in computational chemistry.
The results for the Rabi models (ultrastrong and deep strong coupling) are presented in the following Section~\ref{sec:Results for the Rabi model}.
In Section~\ref{sec:Extensions to higher dimensions}, we discuss the extension of both the koopmon and the bohmion method to higher dimensions.
Finally, we present a brief conclusion and outlook in Section~\ref{sec:Conclusions and Outlook}, while Appendix~\ref{app:sec:Nonadiabatic coupling} expands on certain concepts commonly used in the quantum chemistry of nonadiabatic processes.

\section{The mixed quantum-classical \emph{koopmon} method}\label{sec:The mixed quantum-classical koopmon method}
\subsection{From Koopman wavefunctions to MQC dynamics}\label{sec:QCWE}
Following Sudarshan's work from 1976 \cite{Sudarshan}, the new MQC model \eqref{HybEq1}-\eqref{HybEq3} was obtained in two main steps hinging on the use of Koopman wavefunctions in classical mechanics \cite{Koopman}:
\begin{enumerate}
	\item write the MQC dynamics on the tensor-product space of Koopman and Schr\"odinger wavefunctions; 
	\item make classical phases \emph{unobservable} \cite{Ghose}.
\end{enumerate}
The difficulty in following Sudarshan's plan lies in two crucial points:
i) the dynamics on the tensor-product Hilbert space is far from obvious;
ii) it is not clear what is mathematically meant by `unobservable'.
Indeed, Sudarshan's original model led to interpretative difficulties which attracted various criticism \cite{PeTe}.
In \cite{GBTr22,GBTr21} we circumvented these difficulties by first blending Koopman wavefunctions with van Hove's unitary representations in symplectic geometry \cite{VanHove}, and then by making classical phases an overall symmetry for the resulting dynamics.
 
The first step in this derivation partly follows \cite{JaSu10} and amounts to modifying KvN dynamics in such a way that the Koopman phase is compatible with Hamilton--Jacobi theory.
The resulting equation $i\hbar\partial_t\chi=\{i\hbar H,\chi\}-(p\partial_pH-H)\chi$, which was dubbed \emph{Koopman-van Hove equation} (KvH) \cite{BoGBTr19,TrJo21}, goes back to the work of Bertram Kostant \cite{Kostant}.
In order to tackle Sudarshan's step 1) above, we started with the KvH construction for two classical degrees of freedom and then quantized one of those.
The resulting linear \emph{quantum-classical wave equation} (QCWE) $i\hbar\partial_t\Upsilon=\{i\hbar\widehat{H},\Upsilon\}-(p\partial_p\widehat{H}-\widehat{H})\Upsilon$ for the MQC wavefunction $\Upsilon(q,p,x)$ has a canonical Hamiltonian structure and was largely studied in \cite{GBTr20}.
While it retains Heisenberg's principle, a proof that the QCWE also preserves the sign in the classical sector has so far only been provided for pure-dephasing Hamiltonians, for which a numerical study was presented in \cite{Manfredi23}.
We note that, unlike the classical KvH equation above, the QCWE contains the operator-valued function $\widehat{H}(\bz,\hat q,\hat p)$ (hybrid Hamiltonian operator), which depends on the classical phase-space coordinates as well as the quantum position- and momentum operators $\hat q$ and $\hat p$, respectively.
Later, the model was upgraded \cite{GBTr21} by noticing that Sudarshan's step 2) could be taken by modifying the underlying action principle $\delta\int_{t_1}^{t_2}\int\langle\Upsilon,i\hbar\partial_t\Upsilon-\{i\hbar\widehat{H},\Upsilon\}+(p\partial_p\widehat{H}-\widehat{H})\Upsilon\rangle\,\de q\de p=0$ in such a way that the latter becomes symmetric under phase transformations $\Upsilon\mapsto e^{i\varphi(q,p,t)}\Upsilon$.
Here, $\langle\,,\rangle:=\operatorname{Re}\langle\,|\,\rangle$, where $\langle\,|\,\rangle$ is the Hermitian inner product on the quantum Hilbert space $\mathscr{H}$, with norm $\|\,\|$.
We emphasize that the aforementioned phase symmetry principle is necessary to ensure that the classical Liouville density $f(q,p)$ is positive at all times, so that $f=\|\Upsilon\|^2$.
The application of this phase symmetry requires several steps, which start by isolating the classical phase within the MQC wavefunction $\Upsilon$.
This can be achieved by resorting to conditional wavefunctions, that is by writing $\Upsilon(\bz,x,t)=\chi(\bz,t)\psi(x,t;\bz)$ \cite{Abedi10}.
At each point in time, ${\chi\in L^2(T^*Q)}$ is a Koopman wavefunction, while $\psi$ is a conditional Schr\"odinger wavefunction parameterized by $\bz$, so that $\|\psi\|^2=\int|\psi(\bx,t;\bz)|^2\,\de x=1$ at all times $t$ and for all points $\bz=(q,p)$.
In the remainder of this paper, we will restrict to consider a finite dimensional Hilbert space $\mathscr{H}$ in the quantum sector to avoid dealing with important issues in functional analysis.

Once the classical phase $S(\bz,t)$ is identified by the polar form $\chi=\sqrt{f}e^{iS/\hbar}$, applying the symmetry principle amounts to eliminating this phase from the action principle $\delta\int_{t_1}^{t_2}\!\int\!f(\partial_t S-{\langle\psi,i\hbar\partial_t\psi-\{i\hbar\widehat{H},\psi\}+(p\partial_p\widehat{H}-\widehat{H})\psi+\{S,{\widehat{H}}\}\psi\rangle})\,\de ^2z=0$.
Benefiting from the properties of the \emph{Berry connection $\langle\psi|-i\hbar\nabla\psi\rangle$} \cite{Berry1984}, this step requires introducing diffeomorphic Lagrangian paths $\boldsymbol{\eta}_t\in\operatorname{Diff}(T^*Q)$ on phase space, which are generated by the vector field $\boldsymbol{\cal X}$ via the relation $\dot{\boldsymbol{\eta}_t}=\boldsymbol{\cal X}\circ\boldsymbol{\eta}_t$.
Here, $\operatorname{Diff}(T^*Q)$ denotes the group of diffeomorphism on the classical phase space.
In particular, one combines integration by parts with the continuity equation $\partial_t f+\operatorname{div}(f\boldsymbol{\cal X})=0$ to write $\delta\int_{t_1}^{t_2}\!\int\!f\partial_t S\,\de^2z\de t=-\delta\int_{t_1}^{t_2}\!\int\!f\nabla S\cdot\boldsymbol{\cal X}\,\de^2z\de t$, where $f$ and $\boldsymbol{\cal X}$ are now varied as in Euler--Poincar\'e theory \cite{HoMaRa98}, that is $\delta f=-\operatorname{div}(f\boldsymbol{\cal Y})$ and $\delta\boldsymbol{\cal X}=\partial_t\boldsymbol{\cal Y}+\boldsymbol{\cal X}\cdot\nabla\boldsymbol{\cal Y}-\boldsymbol{\cal Y}\cdot\nabla\boldsymbol{\cal X}$, where $\boldsymbol{\cal Y}=\delta\boldsymbol{\eta}_t\circ\boldsymbol{\eta}_t^{-1}$ is arbitrary.
Finally, one last convenient step consists in writing $\psi(t)=(U_t\psi_0)\circ\boldsymbol\eta_t^{-1}$, without loss of generality \cite{GBTr22}.
Here, $U_t=U_t(\bz)$ is a unitary operator on the quantum Hilbert space that is parameterized by phase-space coordinates.
This last step simply expresses the quantum unitary dynamics in the frame of Lagrangian classical paths.
In this way, the action principle is entirely expressed in terms of $f$ and $\widehat{\cal P}=\Upsilon\Upsilon^\dagger=f\psi\psi^\dagger$ as
\beq\label{VarPrin1}
	\delta\!\int_{t_1}^{t_2}\!\!\left(\int\!\big(f\boldsymbol{\cal A}\cdot\boldsymbol{\cal X}+\langle\widehat{\cal P},i\hbar\hat\xi\rangle\big)\de^2z-h(f,\widehat{\cal P})\!\right)\!\de t
	=0,\quad\,\,\,\,\,
	h(f,\widehat{\cal P})
	=\!\int\!\big\langle f{\widehat{H}}+i\hbar\{\widehat{\cal P},\widehat{H}\}\big\rangle\,\de^2z,
\eeq
which is indeed phase-independent.
Here, ${\boldsymbol{\cal A}=(p,0)}$ is the coordinate representation of the canonical one-form ${\cal A}={\boldsymbol{\cal A}\cdot\de\bz}=p\de q$.
Also, $\hat\xi=(\dot{U}_tU_t^\dagger)\circ\boldsymbol\eta_t^{-1}$ is skew-Hermitian, so that $\delta\hat\xi=\partial_t\hat\Lambda+[\hat\Lambda,\hat\xi]+\boldsymbol{\cal X}\cdot\nabla\hat\Lambda-\boldsymbol{\cal Y}\cdot\nabla\hat\xi$ and $\delta\widehat{\cal P}=[\hat\Lambda,\widehat{\cal P}]-\operatorname{div}(\widehat{\cal P}\boldsymbol{\cal Y})$, where $\hat\Lambda=(\delta{U}_tU_t^\dagger)\circ\boldsymbol\eta_t^{-1}$ is skew-Hermitian and arbitrary.
These variations arise by standard Euler--Poincar\'e reduction from Lagrangian to Eulerian variables \cite{HoMaRa98}.
Indeed, Lagrangian trajectories play a crucial role in the variational problem associated to \eqref{VarPrin1}.
In particular, if $\boldsymbol\eta(\bz_0,t)$ is the diffeomorphic Lagrangian path on phase space and $U(\bz,t)$ is a unitary operator, the Eulerian quantities evolve as follows:
\beq\label{LtoE}
	f
	:=\boldsymbol{\eta}_{t*}f_0,\qquad
	\widehat{\cal P}
	:=\boldsymbol{\eta}_{t*}(U_t\widehat{\cal P}_0U_t^\dagger),\qquad
	\boldsymbol{\cal X}
	:=\dot{\boldsymbol\eta}_t\circ\boldsymbol\eta_t^{-1},\qquad
	\hat\xi
	:=\dot{U}_tU_t^\dagger\circ\boldsymbol\eta_t^{-1},
\eeq
where $\boldsymbol{\eta}_{t*}$ denotes the push-forward by the time-dependent map $\boldsymbol{\eta}_t$.
As shown in \cite{GBTr21,TrGB23}, the variational principle \eqref{VarPrin1}, with the expressions above for the Euler--Poincar\'e variations $\delta f$, $\delta\widehat{\cal P}$, $\delta\boldsymbol{\cal X}$, and $\delta\hat\xi$, return the model equation \eqref{HybEq1} along with its trace, that is $\partial_t{f}+\operatorname{div}\!\big(f\boldsymbol{\cal X}\big)=0$, and the accompanying relations \eqref{HybEq2}-\eqref{HybEq3}.
Notice that the total energy is identified by the Hamiltonian functional $h(f,\widehat{\cal P})$ in \eqref{VarPrin1}, which is given by the sum of usual average $\int\!f\langle\widehat{H}\rangle\de^2z$ with an extra term rendering the \emph{quantum backreaction}, that is the dynamical features specifically associated to the statistical correlations between the quantum and the classical systems.
Here, these correlations appear through the inhomogeneities encoded in the integrand term $\langle i\hbar\{\widehat{\cal P},\widehat{H}\}\rangle={f^{-1}\operatorname{Tr}(\bX_{\widehat{H}}\cdot[i\hbar\widehat{\cal P},\nabla\widehat{\cal P}])}/2$ in \eqref{VarPrin1}.
As noted in \cite{GBTr-fluid,DMCTr24}, an analogous term also appears in the theory of spin-orbit coupling and this analogy hinges on the prominent role of the \emph{Mead connection} $i\hbar f^{-2}[\widehat{\cal P},\nabla\widehat{\cal P}]$ \cite{Me92}.
This and related aspects will be developed elsewhere.

\subsection{Regularization and the particle scheme}\label{sec:Regularization and the particle scheme}
Given the challenging nature of the nonlinear PDEs \eqref{HybEq1}-\eqref{HybEq3}, here we approximate their associated dynamics by devising a particle method that is based on a variational regularization.
The treatment follows an analogous construction developed within our team in the context of quantum hydrodynamics \cite{FoHoTr19,HoRaTr21}, and later applied to Born--Oppenheimer dynamics \cite{RaTr20}.

We begin our discussion by observing that the point-particle solution ansatz \eqref{singansatz1} is prevented by the presence of the backreaction term ${\cal B}(f,\widehat{\cal P})=\int\!{f}^{-1}\langle\widehat{\cal P},i\hbar\{\widehat{\cal P},\widehat{H}\}\big\rangle\,\de^2z$ in the Hamiltonian $h(f,\widehat{\cal P})$ of \eqref{VarPrin1}.
Indeed, the inverse $f^{-1}$ does not allow for a train of delta functions in the classical sector.
Following \cite{FoHoTr19}, here we will restore the singular solutions by a variational regularization procedure.
In particular, we consider a smooth and strictly positive convolution kernel $K^{(\alpha)}\in C^2(T^*Q)$ on phase space depending on a positive parameter $\alpha>0$, so that in the limit $\alpha\to 0$ the sequence $(K^{(\alpha)})$ converges to the delta distribution.
In addition, we ask that $K^{(\alpha)}(\bz)=K^{(\alpha)}(-\bz)$.
In the following, we write $K=K^{(\alpha)}$ and implicitly assume the dependency of the kernel on the parameter $\alpha$.
We use the kernel to introduce the regularized quantities
\begin{align*}
	\bar{f}(\bz,t)
	=\int\!K(\bz-\bz')f(\bz',t)\de^2z',\qquad\quad
	\bar{\cal P}(\bz,t)
	=\int\!K(\bz-\bz')\widehat{\cal P}(\bz',t)\de^2z',	
\end{align*}
and construct the following regularized variational principle:
\beq\label{VarPrin1reg}
	\delta\!\int_{t_1}^{t_2}\!\!\left(\int\!\big(f\boldsymbol{\cal A}\cdot\boldsymbol{\cal X}+\langle\widehat{\cal P},i\hbar\hat\xi-\widehat{H}\rangle\big)\de^2z-{\cal B}_\textit{reg}\right)\!\de t
	=0,\quad\,\,\,\,
	{\cal B}_\textit{reg}
	=\!\int\!\bar{f}^{-1}\langle\bar{\cal P},i\hbar\{\bar{\cal P},\widehat{H}\}\big\rangle\,\de^2z,
\eeq
with the same Euler--Poincar\'e variations as above.
Note that regularized functions are denoted by a bar.
At this point, taking the time derivative of the second in \eqref{LtoE} yields
\beq\label{RegPDE}
	\frac{\partial\widehat{\cal P}}{\partial t}+\operatorname{div}({\widehat{\cal P}}\boldsymbol{\cal X})
	=[\hat\xi,{\widehat{\cal P}}].
\eeq
Moreover, taking variations in \eqref{VarPrin1reg} yields
\beq\label{EPeqns}
	\boldsymbol{{\cal X}}
	=\langle\bX_{\widehat{H}}\rangle+\bX_{\overline{{\delta{\cal B}}/{\delta\bar{f}}}}+\,\Big\langle\bX_{\overline{{\delta{\cal B}}/{\delta\bar{\cal P}}}}\Big\rangle,
	\qquad\qquad 
	i\hbar[\hat\xi,\widehat{\cal P}]
	=\bigg[\widehat{H}+\overline{\frac{\delta{\cal B}}{\delta\bar{\cal P}}},\widehat{\cal P}\bigg].
\eeq
Here, we have used the standard notation for functional derivatives as well as the chain-rule relations ${\delta {\cal B}_\textit{reg}}/{\delta f}=\overline{{\delta {\cal B}(\bar{f},\bar{\cal P})}/{\delta\bar{f}}}$ and ${\delta {\cal B}_\textit{reg}}/{\delta\widehat{\cal P}}=\overline{{\delta {\cal B}(\bar{f},\bar{\cal P})}/{\delta\bar{\cal P}}}$, resulting from ${\cal B}_\textit{reg}(f,\widehat{\cal P})={\cal B}(\bar{f},\bar{\cal P})$.
As a result, the singular ansatz \eqref{singansatz1} is now an exact solution of the equation ${\partial_t}f+\operatorname{div}(f\boldsymbol{{\cal X}})=0$, so that $\dot{\bzeta}_a=\boldsymbol{{\cal X}}(\bzeta_a)$ and one can proceed similarly to the discussion in Section~\ref{sec:Ehrenfest} upon writing the equation $\partial_t\,\widehat{\!\mathscr{P}}+\boldsymbol{{\cal X}}\cdot\nabla\,\widehat{\!\mathscr{P}}=[\hat\xi,\,\widehat{\!\mathscr{P}}]$ for the quantity $\,\widehat{\!\mathscr{P}}=\widehat{\cal P}/f$.
This procedure leads again to a set of closed equations for $(\bzeta_a(t),\hat\varrho_a(t))$, where $\hat\varrho_a=\,\widehat{\!\mathscr{P}}(\bzeta_a)$, thereby extending the multi-trajectory Ehrenfest scheme \eqref{MFeqs}.
Given the complexity of the regularized continuum model \eqref{RegPDE}-\eqref{EPeqns}, however, here we follow a different, yet equivalent method.

By adopting a procedure developed in quantum hydrodynamics \cite{FoHoTr19}, we proceed by inserting the following singular solution ansatz in the regularized variational principle \eqref{VarPrin1reg}:
\begin{align*}
	f(\bz,t)
	=\sum_{a=1}^N w_a\delta(\bz-\bzeta_a(t)),\qquad\quad
	\widehat{\cal P}(\bz,t)
	=\sum_{a=1}^N w_a\hat\varrho_a\delta(\bz-\bzeta_a(t)),	
\end{align*}
This ansatz corresponds to a singular solution momentum map first appeared in the geometry of hydrodynamic systems \cite{HoMa2005} and later extended to the case of complex fluids \cite{HoTr2009}.
We refer to the computational particles identified by the phase-space coordinates $\bzeta_a$ as \emph{koopmons}.
Each koopmon carries a positive weight $w_a$ and a quantum state $\hat\varrho_a$.
Based on this ansatz, \eqref{VarPrin1reg} becomes \cite{TrGB23}
\begin{align*}
	\delta\int_{t_1}^{t_2}\sum_aw_a\bigg(p_a\dot{q}_a+\langle\hat\varrho_a,i\hbar\hat\xi_a-{\widehat{H}}_a\rangle-\frac{1}2\sum_{b}w_b\big\langle i\hbar[{\hat\varrho}_a,{\hat\varrho}_b],\widehat{\cal I}_{ab}\big\rangle\bigg)\de t
	=0,
\end{align*}
where ${\widehat{H}}_a={\widehat{H}}(\bzeta_a)$, and we have introduced the notation
\beq\label{koopint}
	\widehat{\cal I}_{ab}
	:=\frac12\int\frac{K_a\{K_b,\widehat{H}\}-K_b\{K_a,\widehat{H}\}}{\sum_c w_c K_c}\,\de^2z,
	\qquad\text{and}\qquad
	K_s(\bz,t)
	:=K(\bz-\bzeta_s(t)).
\eeq
In addition, we have $\dot{\hat\varrho}_a=[\hat\xi_a,{\hat\varrho}_a]$ and $\hat\xi_a(t)=\hat\xi(\bzeta_a(t),t)$, so that $\delta{\hat\varrho}_a=[\hat\Lambda_a,{\hat\varrho}_a]$ and $\delta\hat\xi_a=\partial_t\hat\Lambda_a+[\hat\Lambda_a,\hat\xi_a]$, where $\hat\Lambda_a$ is arbitrary.
These variations, together with arbitrary $\delta q_a$ and $\delta p_a$, yield the \emph{koopmon scheme}
\beq\label{MFeqs2}
	\dot{q}_a
	=w_a^{-1}{\partial_{p_a\!}h},\qquad
	\dot{p}_a
	=-w_a^{-1}{\partial_{q_a\!}h},\qquad
	i\hbar{\dot{\hat\varrho}_a}
	=w_a^{-1}[{\partial_{\hat\varrho_a\!}h},\hat\varrho_a],
\eeq
where
\beq\label{koopmonHam}
	h
	=\sum_aw_a\langle\hat\varrho_a,{\widehat{H}}_a\rangle+\frac12\sum_{a,b}w_aw_b\big\langle i\hbar[{\hat\varrho}_a,{\hat\varrho}_b],\widehat{\cal I}_{ab}\big\rangle.
\eeq
Notice that discarding the last term in the expression above recovers the multi-trajectory Ehrenfest system \eqref{MFeqs}.
Indeed, both the latter and the koopmon model share the same Hamiltonian structure that is identified by the direct-sum Poisson bracket
\begin{align*}
	\{\!\!\{h,k\}\!\!\}\big((\bzeta_a)_{a=1}^N,(\hat\varrho_a)_{a=1}^N\big)
	=\sum_{a=1}^N\frac1{w_a}\left(\frac{\partial h}{\partial q_a}\frac{\partial k}{\partial p_a}-\frac{\partial h}{\partial p_a}\frac{\partial k}{\partial q_a}-\left\langle\hat\varrho_a,\frac{i}\hbar\left[\frac{\partial h}{\partial\hat\varrho_a},\frac{\partial k}{\partial\hat\varrho_a}\right]\right\rangle\right).	
\end{align*}
Also, equations \eqref{MFeqs2}-\eqref{koopmonHam} inherit all the MQC consistency criteria listed in Section~\ref{sec:Ehrenfest}.
Despite these analogies, we observe that the koopmon scheme differs substantially from \eqref{MFeqs} in that each particle is coupled to all others at all points in time.
This coupling among trajectories reflects the quantum-classical correlation effects and indeed vanishes when $\widehat{H}$ is purely classical (that is $\widehat{H}=H\boldsymbol{1}$) or purely quantum (that is $\nabla\widehat{H}=0$).
Also, we observe that this coupling does not take place via additional force terms in a standard set of Newtonian equations.
Instead, the coupling occurs via an integral involving both the position and momentum of all particles, thereby affecting the evolution of both classical variables.

As previously illustrated in the plan of the paper, so far we have reviewed the main differences between Ehrenfest dynamics and our new MQC model in \eqref{HybEq1}-\eqref{HybEq3}.
In particular, after summarizing the formulation of the underlying continuum PDEs, we have proceeded by presenting the corresponding particle closures.
The reminder of this work presents the results obtained by the implementation of the koopmon scheme along, with a comparison with analogous trajectory methods including the multi-trajectory Ehrenfest method.

\section{Implementation and comparison to other methods}\label{sec:Numerical implementation}
In the remainder of this paper we present the results obtained by the koopmon method and compare them with the analogous results arising from other particle methods in both the MQC and quantum description.
Involving different test cases, these results will also be accompanied by a fully quantum solution, as obtained by a suitable Schr\"odinger solver.
This section discusses the numerical implementation in each of the adopted integration schemes.

The koopmon scheme was implemented in MATLAB to probe its performance across various model systems.
Both MATLAB versions 2021b and 2023b were used during the development of this work.
In what follows we provide a detailed summary of all the different input parameters:
\begin{enumerate}
	\item\emph{Hybrid Hamiltonian}: 
		Each test case requires the operator-valued function $\widehat{H}$ introduced above, which models the quantum-classical interaction.
		Here, we restrict ourselves to cases where the hybrid Hamiltonian models the interaction with a two-level quantum system, and the interaction potential depends only on the classical coordinate $q$ (no momentum coupling).
		Specifically, 
		\beq\label{UsedHam}
			\widehat{H}(q,p)
			=H_C(q,p)\boldsymbol{1}+\widehat{H}_I(q),
		\eeq
		where the phase-space function $H_C$ models the classical system, $\boldsymbol{1}$ is the identity matrix, and the interaction term $\widehat{H}_I(q)\in\operatorname{Her}(\mathbb{C}^2)$ is given by a linear combination of the Pauli matrices $\hat\sigma_x,\hat\sigma_y,\hat\sigma_z\in\mathbb{C}^{2\times 2}$.
		For all test cases in this paper, $\widehat{H}$ will be a smooth function and all physical quantities are expressed in atomic units.
	\item \emph{Number of particles}:
		The number $N\ge 1$ of particles determines the dimensionality of the problem.
		In the classical sector, the arrays $\mathbf{q}:=(q_a)_{a=1}^N$ and $\mathbf{p}:=(p_a)_{a=1}^N$ belong to $\mathbb{R}^N$, whereas the quantum sector is described by the complex-valued array $(\hat\varrho_a)_{a=1}^N\in\mathbb{C}^{N\times 2\times 2}$ of density matrices $\hat\varrho_a\in\operatorname{Her}(\mathbb{C}^2)$.
		The weights were fixed to ${w_a=1/N}$ for all values of $a$.
	\item \emph{Width of the kernel function}:
		Recall that the regularization of the koopmon scheme is based on a kernel function $K$ in phase space.
		In our implementation, $K$ is chosen as the product of one-dimensional normalized Gaussians, that is, $K(q,p)=\tilde K(q)\tilde K(p)$ with
		\begin{align}\label{kernel}
			\tilde K(y)
			:=\frac{1}{\alpha\sqrt{\pi}}\,e^{-y^2/\alpha^2},\quad
			y\in\mathbb{R}.
		\end{align}
		In particular, this choice fulfils the previously introduced property $K\to\delta$ as $\alpha\to0$, such that in the latter limit the method recovers the original nonlinear MQC model \eqref{HybEq1}-\eqref{HybEq3}.
		The kernel width $\alpha>0$ will be defined for each test case.
\end{enumerate}
Notice that, upon using a Gaussian kernel function as defined in input 3, the Poisson brackets in the backreaction integral $\widehat{\cal I}_{ab}$ of \eqref{koopint} vanish in the limit $\alpha\to+\infty$, and as a result, the koopmon equations \eqref{MFeqs2}-\eqref{koopmonHam} recover the multi-trajectory Ehrenfest scheme \eqref{MFeqs}.
	Alternative kernel functions, such as the \emph{Student's t-distribution}, also share this property.

\subsection{Koopmon integration scheme}
The koopmon Hamiltonian $h$ in \eqref{koopmonHam} carries the backreaction integral $\widehat{\cal I}_{ab}$ in \eqref{koopint}, whose integrand contains the products $K_a\{K_b,\widehat{H}\}$ and $K_b\{K_a,\widehat{H}\}$ in the numerator.
The presence of the factors $K_a$ and $K_b$, which are Gaussians of width $\alpha$ and phase-space centers $\bzeta_a$ and $\bzeta_b$, suggests the utilization of a grid that depends on $\alpha$ as well as the arrays $\mathbf{q}$ and $\mathbf{p}$.
For the numerical integration of the phase-space integrals $\widehat{\cal I}_{ab}$, the integration domain is truncated to a rectangular box $\Gamma\subset\mathbb{R}^2$, on which the composite trapezoidal rule is applied.
In particular, the composite trapezoidal rule was implemented using the internal {\tt trapz} command in MATLAB.
The truncation box and the spatial grid sizes are selected, as follows:
\begin{enumerate}\setcounter{enumi}{3}
	\item \emph{Truncation box (time dependent)}:
		Due to the rapid decay of the Gaussian kernel function, values of the integrand for arguments that are far away from the minimum or maximum of the classical coordinates $q_a$ and $p_a$ are negligibly small.
		Thus, the rectangular truncation box is chosen so that it dynamically adapts to these minima/maxima in each time step.
		Upon denoting ${(q_{\operatorname{min}},p_{\operatorname{min}}):=(\min_aq_a,\min_ap_a)}$, and ${(q_{\operatorname{max}},p_{\operatorname{max}}):=(\max_aq_a,\max_ap_a)}$, we use the integration box
		\begin{align*}
			\Gamma(t)
			:=\big[q_{\operatorname{min}}(t)-n_q\sigma_K,q_{\operatorname{max}}(t)+n_q\sigma_K\big]\times\big[p_{\operatorname{min}}(t)-n_p\sigma_K,p_{\operatorname{max}}(t)+n_p\sigma_K\big],			
		\end{align*}
		where $n_q,n_p\ge 1$ are positive integers.
		In this paper, we used ${n_p=n_q=2}$.
		For the test cases treated here, a further extension of the box does not result in significant improvements in accuracy.	
	\item \emph{Spatial grid sizes (time independent)}:
		The grid sizes $dq,dp>0$ are related to the kernel width by $dq=\sigma_K/j_q$ and $dp=\sigma_K/j_p$, where $\sigma_K:=\alpha/\sqrt{2}$ and $j_q,j_p\ge 1$ are positive integers.
		In this paper, we used $j_q=j_p=2$.
		Our study shows that, for the test cases treated here, a further reduction of the grid size does not result in significant improvements in accuracy.	
\end{enumerate}
For the numerical integration of the koopmon equations presented in \eqref{MFeqs2}, we use a fourth-order Runge--Kutta method (RK4), for which we introduce the following input parameters:
\begin{enumerate}\setcounter{enumi}{5}
	\item \emph{Time increment}:
		Depending on the model under consideration (Tully or Rabi), we work with different fixed time steps that ensure that the relative error of the time-evolved energy remains within the tolerance $10^{-2}$.
		The exact numerical values will be given later on in each case of study.
	\item \emph{Initial conditions}:
		The initial condition for the classical position of each trajectory in the koopmon method is sampled from a normal distribution $\mathcal{N}(\mu_q,\sigma_q^2)$, where the mean and the standard deviation $(\mu_q,\sigma_q)\in\mathbb{R}\times\mathbb{R}^+$ are given parameters.
		Analogously, the classical momenta are sampled from a normal distribution $\mathcal{N}(\mu_p,\sigma_p^2)$ for a given $\mu_p\in\mathbb{R}$ and $\sigma_p = 1/(2\sigma_q)$.
		In order to achieve fast convergence of the koopmon method with an increasing number $N$ of trajectories, we do not work with the internal MATLAB command for the generation of normally distributed random numbers.
		Instead, the sampling from the above normal distributions is performed using quasi-random numbers based on the two-dimensional Sobol sequence \cite{Niederreiter}.
		In the quantum sector, the initial condition for the density matrix was uniformly chosen as $\hat\varrho_a=\hat\varrho_0$ for all $a=1,\dots,N$, where $\hat\varrho_0\in\operatorname{Her}(\mathbb{C}^2)$ is a given Hermitian matrix.
\end{enumerate}
%

\subsection{The quantum \emph{bohmion} method}
In order to assess the koopmon method in terms of its accuracy and computational costs, we will compare its results to those obtained by both the MQC Ehrenfest model and fully quantum descriptions.
Indeed, coupled-trajectory particle schemes similar to the MQC koopmon method also appear in fully quantum approaches based on quantum hydrodynamics \cite{Agostini,FoHoTr19}.
In this context, one considers the coupled quantum dynamics of nuclei and electrons and describes the nuclear component by Madelung--Bohm trajectories, while treating the electrons in the standard Schr\"odinger picture.
For example, the \emph{bohmion method} \cite{FoHoTr19,HoRaTr21} is based on exactly the same equations as \eqref{MFeqs2}, although the total energy \eqref{koopmonHam} is replaced by
\begin{align*}
	h
	=\sum_aw_a\langle\hat\varrho_a,{\widehat{H}}_{a}\rangle+\frac{\hbar^2}{8M}\sum_{a,b}w_aw_b (2\langle\hat\varrho_{a},\hat\varrho_{b}\rangle-1){\mathscr{ I}}_{ab},	
\end{align*}
where
\beq\label{bohmintegs}
	{\mathscr{I}}_{ab}
	=\int\frac{\partial_r{\sf K}(r-q_a)\partial_r{\sf K}(r-q_b)}{\sum_c w_c{\sf K}(r-q_c)}\,\de r.
\eeq
Also, $M>0$ is the mass of the nuclear subsystem and ${\sf K}$ is a kernel function analogous to $K$, although now defined only on configuration space.
Once again, ${\sf K}$ depends on a parameter $\alpha>0$ such that ${\sf K}\to\delta$ as $\alpha\to0$.
In the latter limit, unlike koopmon dynamics, the bohmion method recovers the exact fully quantum dynamics.
Alternatively, the limit $\alpha\to+\infty$ yields the MQC Ehrenfest equations \eqref{MFeqs}, similarly to the koopmon scheme.
We also emphasize that, like other approaches based on quantum hydrodynamics \cite{Hall}, the bohmion method fails to reduce to standard quantum and classical mechanics in the absence of a potential coupling nuclear and electronic dynamics.
This fact is due to the persistence of the integral ${\mathscr{ I}}_{ab}$ and should not come as a surprise, since the bohmions provide a fully quantum description, not a MQC scheme.

The bohmion method was successfully benchmarked in \cite{HoRaTr21} against well known test cases in nonadiabatic molecular dynamics, thereby demonstrating great accuracy levels even in challenging situations involving avoided crossings and extended coupling regions.
While the integral ${\mathscr{ I}}_{ab}$ is on configuration space, it may still require substantial computational resources compared to the phase-space integral in \eqref{koopint}.
In particular, the computational cost depends on the values of $\alpha$ that are required to achieve sufficient levels of accuracy in each case.
Following the algorithm presented in \cite{HoRaTr21}, we implemented the bohmion method using the (one-dimensional) composite midpoint rule and an RK4 integrator based on the same input parameters outlined above.
Moreover, the chosen kernel ${\sf K}(r)=\tilde K(r)$ is as in \eqref{kernel}.
Since both the koopmon and bohmion methods share the same parameters, our numerical study allows for an effective comparison of both methods.

\subsection{Schr\"odinger solver and output visualization}\label{sub:Quantum reference solver and visualization}
Besides the comparison to the Ehrenfest and bohmion methods, it is important to also compare koopmon dynamics to the solutions of the Schr\"odinger equation for the corresponding fully quantum system.
These quantum solutions are generated using the well-known \emph{split operator Fourier transform (SOFT) method} \cite{Lubich,Greene17}.
The corresponding quantum Hamiltonian is derived by replacing the classical phase-space coordinates $(q,p)$ in the hybrid Hamiltonian $\widehat{H}(q,p)$ with quantum position and momentum operators, respectively.
In our SOFT simulations, the initial wavepacket is chosen as ${\Psi_0:=\psi_0v_0\in L^2(\mathbb{R},\mathbb{C}^2)}$, where $v_0\in\mathbb{C}^2$.
Depending on the specific model, $v_0$ is chosen as the eigenvector of the Hamiltonian $\hat H_I(\mu_q)$ either in the ground or the excited state.
Additionally, the initial wavefunction $\psi_0$ is a normalized Gaussian wavepacket centered at $(\mu_q,\mu_p)\in\mathbb{R}^2$, that is,
$
	\psi_0(r)
	=(\gamma/\pi)^{1/4}\exp\left(\mu_p(r-\mu_q)-{\gamma}(r-\mu_q)^2/2\right)
$,
with $\gamma=(2\sigma_q^2)^{-1}$.
In other words, the initialization for the koopmon method in the classical sector, as described in input 7, can be viewed as a sampling from the Wigner distribution $W_{\psi_0}\in L^1(\mathbb{R}^2)$.
Moreover, in the quantum sector, we have $\hat\rho_0=v_0v_0^\dagger$.

As anticipated, we will analyze koopmon dynamics and compare its features to those appearing in SOFT, Ehrenfest, and bohmion dynamics.
Various quantities are visualized in both the classical and quantum sectors.
For the particle methods, we present the classical Liouville density by plotting the particle clouds $(\mathbf{q}(t),\mathbf{p}(t))\in\mathbb{R}^N\times\mathbb{R}^N$ at different times.
In SOFT simulations, we plot the Wigner distribution of the time-evolved wavepacket $\Psi(t)$.
More precisely, for the vector-valued wavefunction $\Psi=(\psi_1,\psi_2)^T$ we compute $W_\Psi:=W_{\psi_1}+W_{\psi_2}$, where $W_{\psi_j},\,j\in\{1,2\}$, denotes the usual Wigner distribution of the one-dimensional component $\psi_j\in L^2(\mathbb{R})$, given by
\begin{align*}
	W_{\psi_j}(q,p)
	:=\frac{1}{\pi\hbar}\int_{-\infty}^\infty\psi_j^*(q+y)\psi_j(q-y)e^{2ipy/\hbar}\,\mathrm{d}y,
	\qquad\text{for all $(q,p)\in\mathbb{R}^2$.}
\end{align*}
In the reminder of this paper, the term \emph{wavepacket} refers indifferently to either the wavefunction representation in configuration space or its Wigner distribution in phase space.
Notice that, unlike the classical Liouville density, the quantum Wigner distribution can take negative values.
As we will see, areas with a negative sign in the Wigner distribution are particularly crucial for analyzing the koopmon method.
We recall that the negative values developed by Wigner dynamics are distinctive tracks of purely quantum features, typically related to interference effects.
To enhance the visual comparison of the Wigner distribution with the classical particle clouds, the latter will be supplemented with the regularized density
$
	D(\bz,t)
	:=(K^{(\Delta)}*f)(\bz,t)=\sum_aw_aK^{(\Delta)}(\bz-\bzeta_a(t)),
$
where $K^{(\Delta)}$ is a Gaussian kernel whose width $\Delta>0$ is a purely visualization parameter, independent of the model parameter $\alpha$.
Here, we have fixed the value $\Delta=0.25$ for all test cases.
In the quantum sector, we plot the evolution of the quantum purity $\|\hat\rho(t)\|^2=\operatorname{Tr}\hat\rho(t)^2$, which is a basis-independent measure of decoherence.
Note that in SOFT simulations we have $\hat\rho(t)=\int\Psi(r,t)\Psi(r,t)^\dagger\,\de r$, while $\hat\rho(t)=\sum_a w_a\hat\varrho_a(t)$ for the particle methods.

\section{Results for the Tully models}\label{sec:Results for the Tully models}
We start the presentation of our numerical results for the Tully models, originally introduced by John Tully to probe his \emph{fewest switches surface hopping algorithm} \cite{Tully90}.
Over the decades, these models have evolved into essential tools for testing novel MQC algorithms in computational chemistry \cite{Ibele20,Nassimi}.
Each Tully model serves as a platform for simulating different nonadiabatic molecular processes; see Appendix~\ref{app:sec:Nonadiabatic coupling} for further details from a chemical perspective.
For example, the bohmion scheme was successfully tested against all three Tully models, for which it exhibited great accuracy levels with a regularization parameter $\alpha\leq1/20$.
The first two models are designed to simulate the dynamics of a wavepacket entering a coupling region (I: single avoided crossing, II: double avoided crossing) on the lower surface.
By passing through these coupling areas, parts of the wavepacket move from the lower to the upper surface.
The third model is known to present a typically challenging test case, modeling not only the transmission of wavepackets as in Tully~I and II, but also involving reflection effects.
All results presented in this paper are based on initializations in low (I+II) or intermediate (III) momentum regimes, which are known to be more challenging test cases than those with high initial momentum.

\begin{figure}[h]
	\centering
	\includegraphics[width=\textwidth]{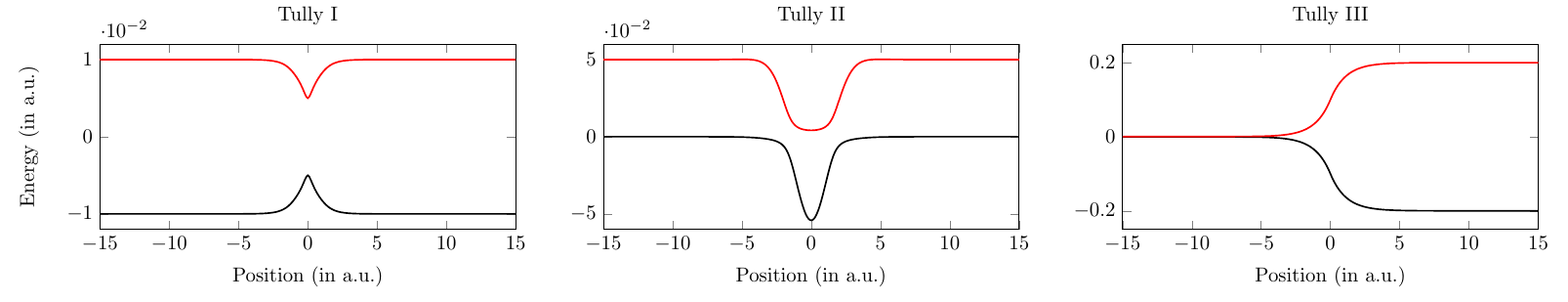}\vspace{-.5cm}
	\caption{\footnotesize
		Potential energy surfaces for the Tully models.
	}\label{pes_T}
\end{figure}

In the MQC implementation, each Tully model corresponds to a different Hamiltonian of the type \eqref{UsedHam}.
In particular,
\beq\label{TullyHam}
	\widehat{H}(q,p)
	=\left(\frac{p^2}{2M}+H_0(q)\right)\!\boldsymbol{1}+H_1(q)\hat\sigma_x+H_3(q)\hat\sigma_z,
\eeq
with model-dependent functions $H_0$, $H_1$, and $H_3$.
Otherwise, except for the initial conditions, the same input parameters listed in Section~\ref{sec:Numerical implementation} were used across all three Tully models.
These are:
Input 1) The classical subsystem is characterized by an effective mass $M=2000$.
Input 2) The number of particles was set to $N=1000$.
Input 3) The kernel width was set to $\alpha=0.325$.
Input 6) The time increment was set to $dt=2$.
Input 7) The width of the initial Gaussian wavepacket is linked to the initial momentum via the relation $\sigma_q=20/(\sqrt{2}\mu_p)$ \cite{Agostini,HoRaTr21}.
The parameter values in inputs 2) and 6) allow us to run the simulations on a laptop machine, thereby avoiding the need for computer clusters.
As for the modeling parameter $\alpha$, values as high as $\alpha=0.75$ were also used on some occasions and we found that $\alpha=0.5$ already allows for a good comparison with the fully quantum results in all the cases considered in this work.
The same holds by taking the number of particles down to $N=500$.
The final values $N=1000$ and $\alpha=0.325$, chosen to be the same for all Tully models, are the result of a compromise between the computational cost associated to the evaluation of the integrals ${\cal I}_{ab}$ in \eqref{koopint} and the need for satisfactory accuracy levels for all Tully models.
Notice that the chosen value of $\alpha$ is well above the $\alpha=1/20$ used in previous bohmion simulations of the Tully models \cite{HoRaTr21}.
Indeed, in that case, such low values of $\alpha$ and their corresponding fine integration grid are required by the complexity of the integrals ${\mathscr{I}}_{ab}$ in \eqref{bohmintegs}, which comprise the regularization of the quantum potential.

\begin{figure}[t]
	\begin{center}
		\addtolength{\leftskip}{-20mm}
		\addtolength{\rightskip}{-20mm}
		\includegraphics[width=172mm]{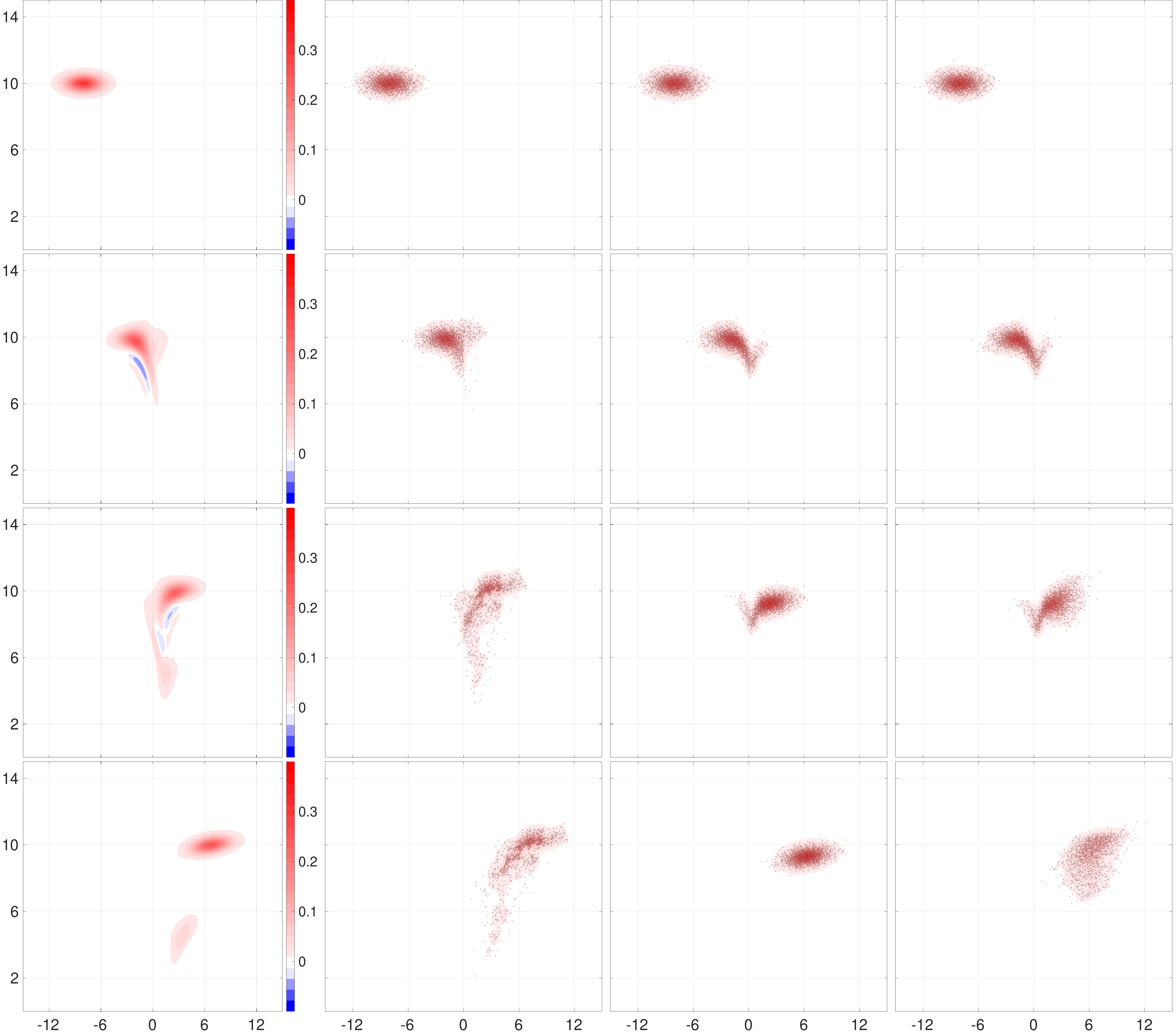}\vspace{-.25cm}
		\caption{\footnotesize
			Time evolution in the classical sector for Tully~I.
			Each column (1st: SOFT, 2nd: koopmons, 3rd: Ehrenfest, 4th: bohmions) shows a snapshot of the distribution in phase space.
			Momentum on the vertical axis, position on the horizontal axis (both in atomic units).
			The rows correspond to times $t=0$, $1280$, $2130$, and $3000$ (top to bottom, all in atomic units).
			In addition, $N=1000$ and $\alpha=0.325$.\vspace{-5.5mm}
		}\label{T1_c}
	\end{center}
\end{figure}

\subsection{Tully~I: single avoided crossing}
The hybrid Hamiltonian for Tully~I results from \eqref{TullyHam} upon letting
\begin{align*}
	H_0(q)
	=0,\qquad
	H_1(q)
	=ce^{-dq^2},\qquad
	H_3(q)
	=a\operatorname{sgn}(q)\big(1-e^{-b|q|}\big),
\end{align*}
with $a=0.01,\,b=1.6,\,c=0.005,\,d=1$.
A plot of the corresponding potential energy surfaces (PESs) can be found in the left box of Figure~\ref{pes_T}.
These PESs are defined as the $q$-dependent eigenvalues of the $2\times2$ matrix ${H_0(q){\boldsymbol{1}}+H_1(q)\hat\sigma_x+H_3(q)\hat\sigma_z}$.
In the classical sector, the initial mean values for position and momentum (Input 7) are chosen as $(\mu_q,\mu_p)=(-8,10)$.
Note that the chosen parameters result in the initial kinetic energy $E_{\operatorname{kin}}(0)=p_0^2/(2M)=0.025$, which lies above the upper PES.
Moreover, in the quantum sector, the initial density matrix represents the ground state, which corresponds to the choice $\hat\rho_0=v_0v_0^\dagger$ for $v_0=(1,0)^T$.
The same set of model parameters has been used previously by other authors \cite{Agostini,HoRaTr21}.

\paragraph{Classical sector.}
Each row of Figure~\ref{T1_c} shows a snapshot of the phase-space distribution obtained by a different method, up to the final time $t=3000$ of the simulation.
As described in Section~\ref{sub:Quantum reference solver and visualization}, the particle clouds are overlapped to the convoluted density ${D:=K^{(\Delta)}*f}$ for visualization purposes.
As we can see for the Wigner distribution in the left column, the wavepacket starts from its initial position and evolves towards the right due to its positive initial momentum.
When reaching the avoided crossing, nonadiabatic transitions occur, causing the wavepacket to split.
Consequently, two regions with different momenta emerge, which are seen in the panel for the final time.
We recall that energy crossings correspond to spectral degeneracies typically manifesting as points of intersection of the PESs.
The regions where the PESs become close are referred to as \emph{avoided crossings}; see Figure~\ref{pes_T}.

For early simulation times, all particle methods are able to reproduce the shape of the Wigner distribution.
Afterwards, the nonadiabatic transitions make it more challenging for all particle methods to accurately capture the dynamics.
Also, we note that the Wigner distribution starts to develop negative values, identifying purely quantum effects.
While Ehrenfest dynamics (third column) fail to reproduce the final configuration, we can see that for the bohmion method (last column) some particles have shifted towards regions of lower momentum.
However, as hinted above, the chosen model parameter $\alpha=0.325$ is not small enough for bohmions to allow for the splitting to be captured precisely.
The koopmon method (second column) yields the best results at all times.
In particular, as we can see in the snapshot for the final time (last row), the number of particles that have moved into the region of lower momenta (below $p\approx 6$) is large enough to show the onset of a splitting, which, however, never becomes distinct.
We have seen similar behaviour in other test cases where no significant mass transfer is observed in the exact quantum dynamics.
Having said that, we also point out that other investigations within the intermediate momentum regimes show that, whenever the particles are distributed more equally across the two PES, the splitting is more distinct and the koopmons succeed in achieving more accurate results.

\begin{figure}[h]
	\centering
	\includegraphics[width=\textwidth]{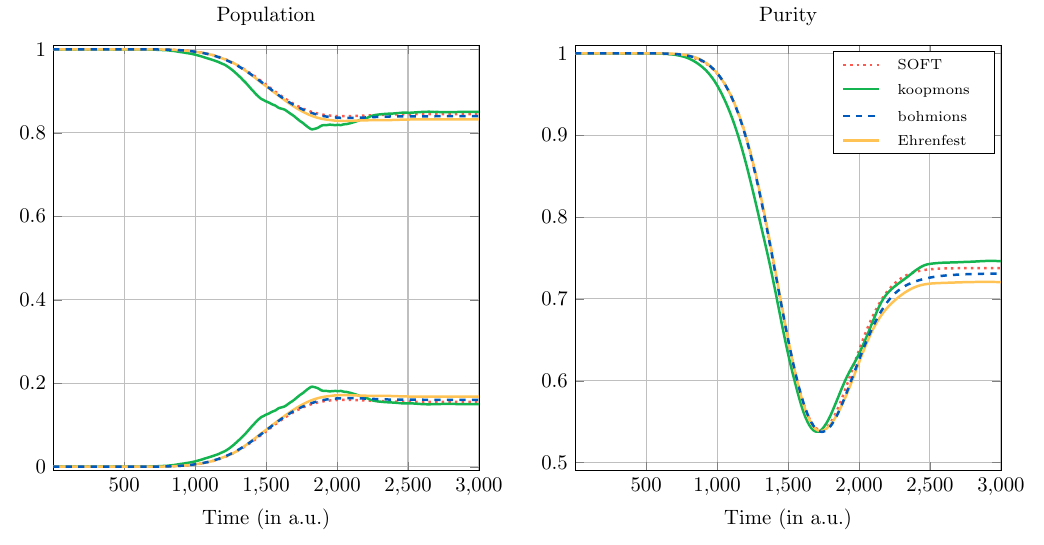}\vspace{-.25cm}
	\caption{\footnotesize
		Population (left) and quantum purity (right) for Tully~I.
	}\label{T1_pp}
\end{figure}

\paragraph{Quantum sector.}
The results in the quantum sector can be found in Figure~\ref{T1_pp}.
In addition to the purity plots at the right-hand side, we also illustrate the \emph{populations}, which represent the proportion of the system in either the ground or the excited state.
In the SOFT simulations, the population on the ground state is calculated as $P_1(t)=\int|\langle v_1(r)|\Psi(r,t)\rangle|^2\,\mathrm{d}r$, where $v_1(r)\in\mathbb{C}^2$ denotes the eigenvector of $\widehat H_I(r)$ corresponding to the smaller eigenvalue.
For the particle methods, we have $P_1(t)=\sum_a w_a\langle v_{a,1}(t)|\hat\varrho_a(t)v_{a,1}(t)\rangle$, where $v_{a,1}(t):=v_1(q_a(t))$.
In other words, the population of the ground state is given by the weighted sum of the individual population levels $P_{1,a}(t):=\langle v_{a,1}(t)|\hat\varrho_a(t)v_{a,1}(t)\rangle$ of each particle.
The population on the excited state is given by $P_2(t)=1-P_1(t)$.

The population plots on the left show that between times $t\approx 1000$ and $t\approx 2250$, the koopmon method (green) has difficulties to accurately reproduce the values obtained by the quantum solution (red).
We observe that in this time region the negative values of the Wigner function become distinctively visible in the fully quantum simulation: these negative values produce purely quantum effects that are not captured by the MQC koopmon method.
However, for the final times, similar to the Ehrenfest (orange) and bohmion (blue) simulations, the koopmons successfully capture the correct population levels.
The purity results show that the koopmon method yields good agreement with the quantum solution at all times.
In particular, at the final time, we can see that the error for both the bohmions and the koopmons is very similar.
However, the koopmons slightly overestimate, while the bohmions underestimate the correct purity levels.
In comparison, Ehrenfest dynamics exhibit small overdecoherence, indicated by lower purity values.

\subsection{Tully~II: dual avoided crossing}
As we will see, the phenomenology accompanying the dual avoided crossing makes Tully II the most challenging Tully problem for the koopmons.
In this case, the hybrid Hamiltonian results from \eqref{TullyHam} upon letting
\begin{align*}
	H_0(q)=-H_3(q)
	=e_0-ae^{-bq^2},\qquad
	H_1(q)
	=c e^{-dq^2},
\end{align*}
with $a=0.05,\,b=0.28,\,c=0.015,\,d=0.06,\,e_0=0.025$.
A plot of the corresponding PESs can be found in the middle panel of Figure~\ref{pes_T}.
In the classical sector, the initial mean values for position and momentum (Input 7) are chosen as $(\mu_q,\mu_p)=(-8,16)$.
Note that the chosen parameters result in the initial kinetic energy $E_{\operatorname{kin}}(0)=0.064$, which again lies above the upper PES; see Figure~\ref{pes_T}.
Moreover, in the quantum sector, the initial density matrix represents the ground state, which corresponds to $\hat\rho_0=v_0v_0^\dagger$ for $v_0=(1,0)^T$.

\paragraph{Classical sector.}
\begin{figure}[h]
	\begin{center}
		\addtolength{\leftskip}{-20mm}
		\addtolength{\rightskip}{-20mm}
		\includegraphics[width=172mm]{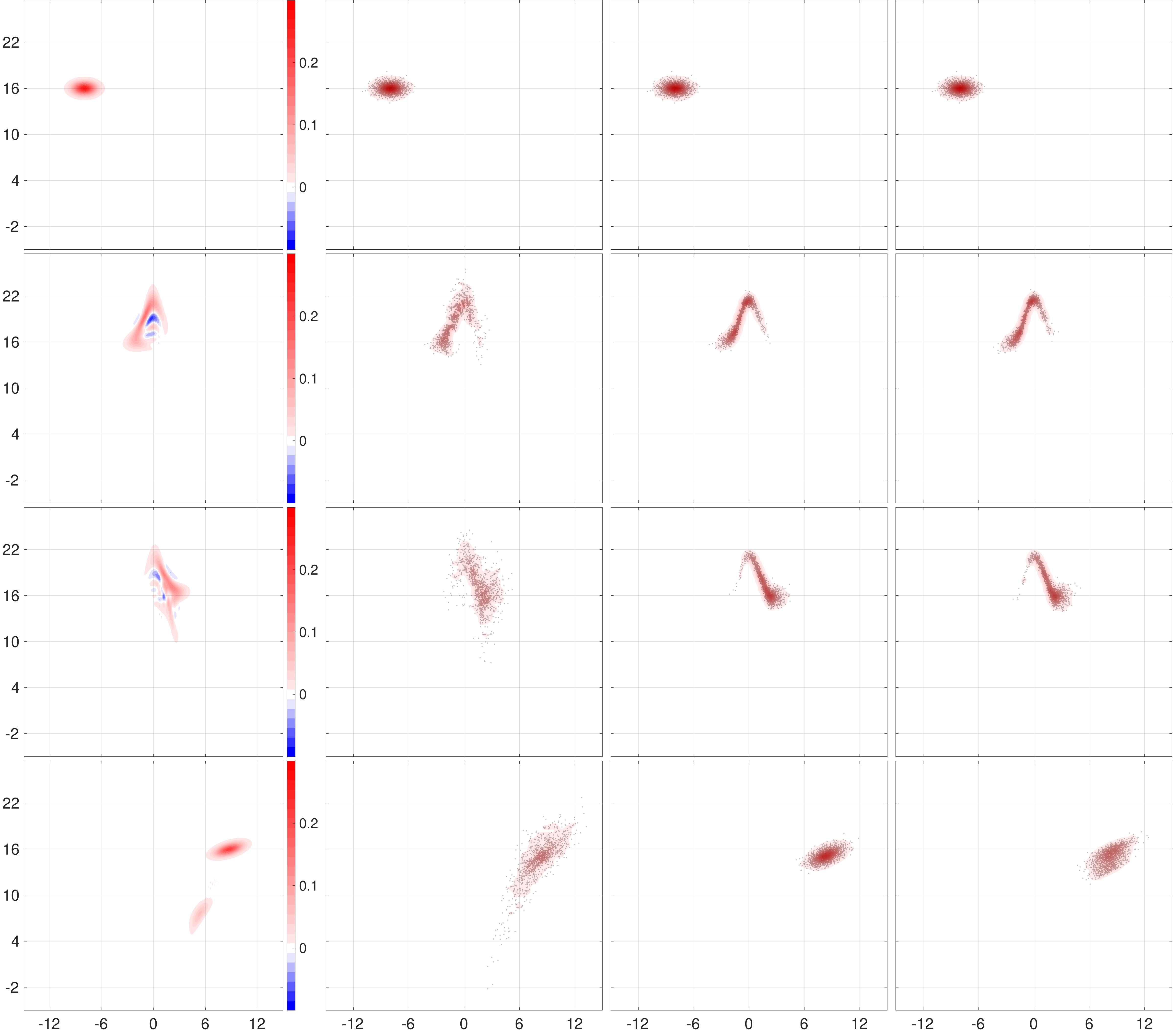}\vspace{-.25cm}
		\caption{\footnotesize
			Time evolution in the classical sector for Tully~II.
			Each column (1st: SOFT, 2nd: koopmons, 3rd: Ehrenfest, 4th: bohmions) shows a snapshot of the distribution in phase space.
			Momentum on the vertical axis, position on the horizontal axis (both in atomic units).
			The rows correspond to times $t=0$, $860$, $1140$, and $2000$ (top to bottom, all in atomic units).
			In addition, $N=1000$ and $\alpha=0.325$.\vspace{-5.5mm}
		}\label{T2_c}
	\end{center}
\end{figure}
Figure~\ref{T2_c} shows the snapshots up to the final time $t=2000$ of the simulation.
At first glance, the dynamics of the quantum solution appears to share qualitative similarities with the results obtained for Tully~I.
However, it is crucial to emphasize that, as seen in Figure~\ref{pes_T}, the PESs of Tully~II comprises two avoided crossings, located on opposite sides of the origin.
Unlike Tully~I, the second avoided crossing in Tully~II introduces additional transition effects as the wavepacket travels through the crossing area.
Furthermore, there is now the possibility of wavepackets recombining.
The strength of these effects depends on the momentum of the initial wavepacket.
In the case of low momentum considered here, this model presents a more challenging test case as compared to Tully~I.

The second row shows that all particle methods yield similar results, accurately reproducing areas where the Wigner distribution attains positive values.
Noticeable differences in the density of the particle clouds can be seen between Ehrenfest or bohmion dynamics, where particles are denser, and the koopmon method, for which the particles show a wider spread.

Similar to Tully~I, both Ehrenfest and bohmion simulations fall short in reproducing the final wavepacket configuration.
In contrast, for the koopmon method some of the particles are observed to move towards regions of lower momentum, resulting in a spreading of the particles rather than a definite splitting.
We would like to note that further simulations indicate that increasing the number of particles in the koopmon method does not alleviate this effect.
Indeed, the particle increase may need to be accompanied by a lower value of $\alpha$.
We emphasize that, similarly to what happens in Tully I, other investigations within the intermediate momentum regimes show that, whenever the particles are distributed more equally across the two PES, the splitting is more distinct and the koopmons succeed in achieving more accurate results.

\begin{figure}[h]
	\centering
	\includegraphics[width=\textwidth]{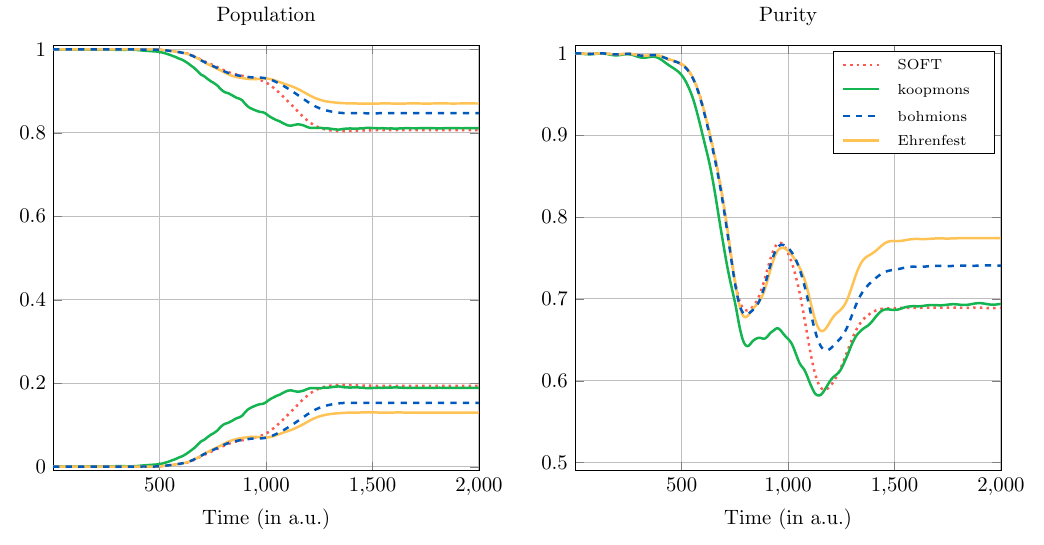}\vspace{-.25cm}
	\caption{\footnotesize
		Population (left) and quantum purity (right) for the model system Tully~II.
	}\label{T2_pp}
\end{figure}

\paragraph{Quantum sector.}
The results in the quantum sector can be found in Figure~\ref{T2_pp}.
Similar to the results for Tully~I, the population plots indicate that the koopmon method has difficulties in accurately reproducing the correct population levels for intermediate times, starting at $t\approx 500$.
At the same time, a slight discrepancy with the fully quantum solution also develops in the purity levels, and this discrepancy gets bigger around $t\approx 750$ while it considerably decreases after $t\approx 1200$.
Once again, these intermediate times between 500 and 1200 correspond to those in which the Wigner function of the fully quantum simulation develops noticeable negative values, thereby leading to quantum effects that are not captured by koopmon dynamics.

A somewhat similar mismatch is also observed in the Ehrenfest simulation, although in this case the discrepancy develops at later times (around $t\approx 1000$), well after the negative phase-space regions have formed, and never decreases even after these negative regions become inappreciable.
The bohmion simulation yields results that are not too dissimilar from Ehrenfest dynamics, although in this case the presence of quantum potential contributions takes the population and purity levels closer to the fully quantum solution.
As shown in \cite{HoRaTr21}, better bohmion performance requires much smaller values of $\alpha$, which are accompanied by higher computational costs.
Unlike bohmions and Ehrenfest, the koopmon method provides final population and purity levels that are in good agreement with the quantum solution.
Although the koopmon method deviates from the quantum solution at intermediate times, it reproduces the correct values after purity has reached its minimum.

\subsection{Tully~III: extended coupling region}
While Tully~III is probably known as the most challenging Tully model for conventional approaches, in this case the koopmon method is particularly successful.
The hybrid Hamiltonian for Tully~III results from \eqref{TullyHam} upon letting
\begin{align*}
	H_0(q)
	=0,\qquad
	H_1(q)
	=
	\begin{cases}
		b(2-e^{-cq}) &\text{if $q>0$},\\
		be^{cq} &\text{if $q\le 0$},
	\end{cases}
	\qquad
	H_3(q)
	=a,
\end{align*}
with $a=0.0006,\,b=0.1,\,c=0.9$.
A plot of the corresponding PESs can be found in the third box of Figure~\ref{pes_T}.
Starting in the ground state, the initial wavepacket moves towards the extended coupling region, causing parts of the wavepacket to populate the upper surface.
However, due to insufficient momentum, the part on the upper surface lacks the energy required to surpass the potential barrier, resulting in a reflection.
This reflection causes the wavepacket to travel through the extended coupling region again, producing further nonadiabatic transitions.
In the classical sector, the initial mean values for position and momentum (Input 7) are chosen as $(\mu_q,\mu_p)=(-15,20)$, following \cite{Subotnik}.
Note that the chosen parameters result in the initial kinetic energy $E_{\operatorname{kin}}(0)=0.1$, which lies between the minimum and the maximum values of the upper PES; see Figure~\ref{pes_T}.
The intermediate momentum value is chosen to allow for a visual comparison with the fully quantum Wigner dynamics over a sufficiently long simulation time, which would be more difficult to plot for lower momentum values.
Moreover, in the quantum sector, the initial density matrix represents the ground state, which in this case corresponds to the choice $\hat\rho_0=v_0v_0^\dagger$ for $v_0=(0,1)^T$.
\begin{figure}[h]
	\begin{center}
		\addtolength{\leftskip}{-20mm}
		\addtolength{\rightskip}{-20mm}
		\includegraphics[width=172mm]{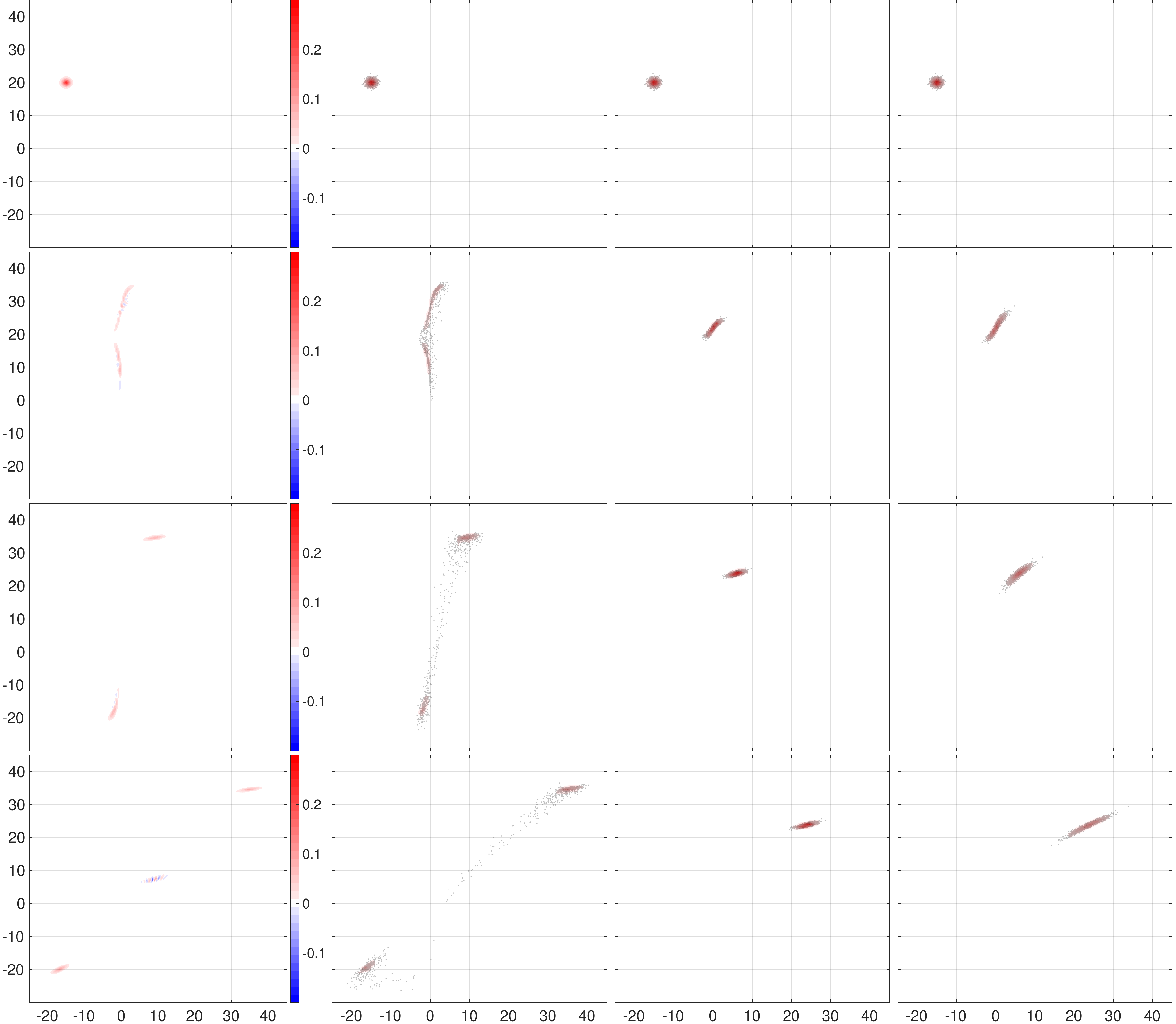}\vspace{-.25cm}
		\caption{\footnotesize
			Time evolution in the classical sector for Tully~III.
			Each column (1st: SOFT, 2nd: koopmons, 3rd: Ehrenfest, 4th: bohmions) shows a snapshot of the distribution in phase space.
			Momentum on the vertical axis, position on the horizontal axis (both in atomic units).
			The rows correspond to times $t=0$, $1500$, $2000$, and $3500$ (top to bottom, all in atomic units).
			In addition, $N=1000$ and $\alpha=0.325$.\vspace{-5.5mm}
		}\label{T3_c}
	\end{center}
\end{figure}

\paragraph{Classical sector.}
Figure~\ref{T3_c} shows the snapshots up to the time $t=3500$.
At time $t=1500$ (second row), the Wigner distribution shows a branching.
The koopmon method captures this branching very well, while the other particle methods only show a slight squeezing of their particle clouds.
At later times, two distinct parts of the wavepacket emerge -- one with positive momentum at the top and one with negative momentum at the bottom.
Note that the negative momentum corresponds to the reflected part of the wavepacket.
The koopmon method reproduces this formation very well.
Also, in the particle cloud for the last snapshot, the lower part is well separated from the rest, while the upper part appears connected to a trail of particles extending toward the origin of phase space.
This is not surprising, since in the last snapshot we can see that the Wigner distribution also forms an alternation of positive and negative values around the phase-space point (10,10).
In contrast to the koopmon dynamics, both Ehrenfest and bohmion simulations only show an increase in the squeezing of the wavepacket over time, without capturing the splitting.

\begin{figure}[h]
	\centering
	\includegraphics[width=\textwidth]{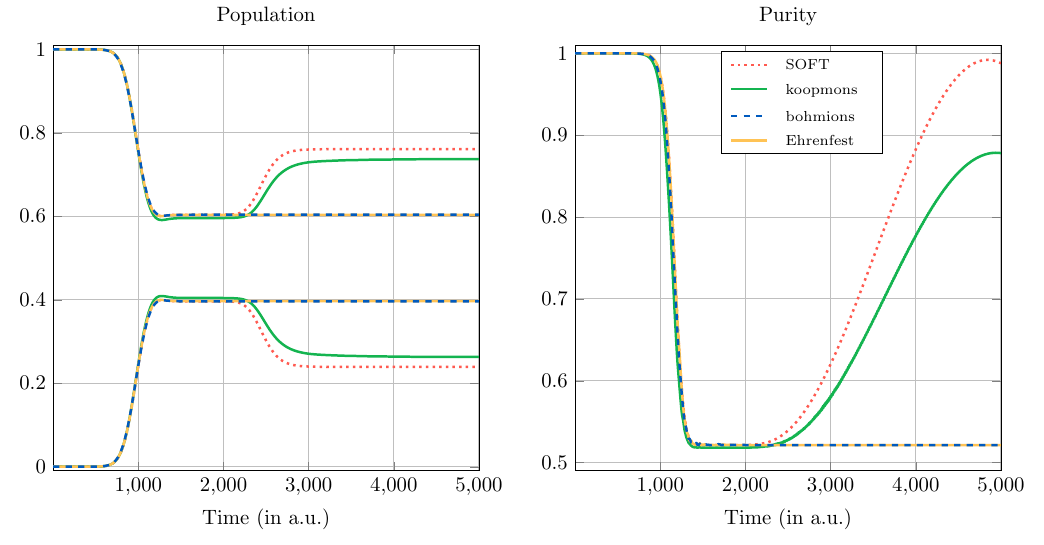}\vspace{-.25cm}
	\caption{\footnotesize
		Population (left) and quantum purity (right) for the model system Tully~III.
	}\label{T3_pp}
\end{figure}

\paragraph{Quantum sector.}
Figure~\ref{T3_pp} shows the superior performance of the koopmon method over the other particle methods in the quantum sector.
While all methods accurately reproduce the population levels up to around $t\approx 2000$, only the koopmon method maintains accurate values at later times.
In the Ehrenfest and bohmion simulations, the populations remain constant for later times.
Even though the koopmon method does not exactly match the levels of the quantum solution, it stands out as it shows an increase in population (for the ground state) around $t\approx 2200$, which is a result of the reflection of the wavepacket on the upper surface.

A similar pattern is observed in the purity plots on the right-hand side.
Until approximately $t\approx 2000$, the curves of all methods overlap.
However, after this point in time, the purity levels of the fully quantum solution go up again -- a pattern captured exclusively by the koopmon method.
Although the koopmon method cannot exactly reproduce the levels of the quantum solution, it qualitatively captures the dynamics.
In contrast, the results for Ehrenfest and bohmion dynamics remain constant for later times.
Once again, we observe that the Wigner distribution of the fully quantum solution presents negative values already at $t\approx 1500$.
Unlike the previous study cases, these negative values do not disappear over time, but instead persist for the entire duration of the simulation and become particularly appreciable after $t\approx 3000$, when the purity levels predicted by the koopmon method start developing a higher discrepancy from the corresponding fully quantum profile.
In this case, the discrepancy may be reduced by decreasing the value of $\alpha$, although at the expense of increasing the computational costs.

\section{Results for the Rabi model}\label{sec:Results for the Rabi model}
\begin{figure}[h]
	\centering
	\includegraphics[scale=.85]{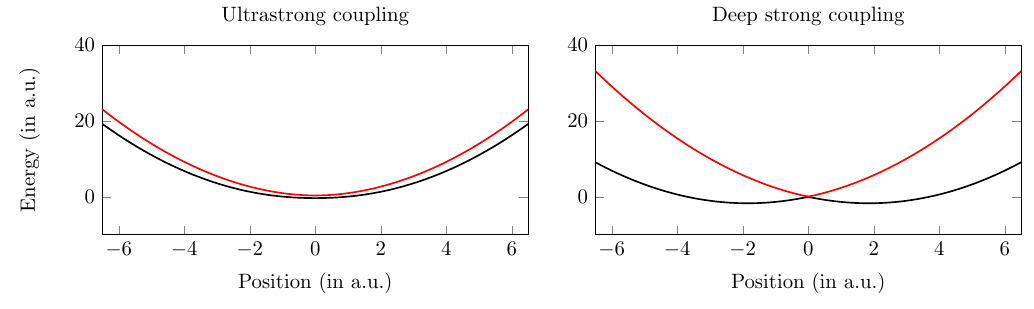}\vspace{-.25cm}
	\caption{\footnotesize
		Potential energy surfaces for the Rabi models.
	}\label{pes_R}
\end{figure}

The two test cases in this section are taken from the physics literature and involve the interaction between a classical harmonic oscillator and a two-level quantum system.
In the fully quantum domain, these models are commonly referred to as \emph{quantum Rabi models} \cite{Rabi}.
These models arise in a variety of contexts, from quantum optics to condensed matter physics, and the theory of open quantum systems; we refer to \cite{PeEtAl15} for a discussion of the various parameter regimes and their relevance.
A quantum-classical study of its pure-dephasing limits was recently presented in \cite{Manfredi23} by resorting to the quantum-classical wave equation discussed in Section~\ref{sec:QCWE}.
Our study covers the ultrastrong and deep strong coupling regimes whose potential energy surfaces are presented in Figure~\ref{pes_R}.
As we will see, MQC treatments are hardly applicable in these regimes.
Nevertheless, we will show that the koopmon method succeeds in reproducing parts of the fully quantum solution.

As in the case of the Tully models presented above, some of the input parameters listed in Section~\ref{sec:Numerical implementation} were used for both Rabi models.
These are as follows:
The classical subsystem is described by a Hamiltonian of the type \eqref{UsedHam}, with $H_C(q,p)=(p^2/M+M\omega^2q^2)/2$ with $M=\omega=1$.
Moreover, 
\beq\label{RabIntHam}
	\widehat{H}_I(q)
	=\gamma q\hat\sigma_z+C_0\hat\sigma_x
\eeq
with model-dependent coupling parameter $\gamma>0$ and $C_0>0$.
The latter represents an external magnetic field driving the spin variable.
The number of particles was set to $N=500$, while the kernel width was set to $\alpha=0.5$.
These parameter values were chosen here to reduce the numerical costs while retaining sufficient levels of accuracy.
We emphasize that the flexibility in these parameters offers a tailored approach based on the characteristics of the dynamics in each model.
In addition, the time increment was set to $dt=0.05$ (Input 6) and the width of the initial Gaussian wavepacket was taken as $\sigma_q=1/\sqrt{2}$.

\subsection{Ultrastrong coupling regime}
We consider a case of ultrastrong coupling regime by setting the constant $\gamma=0.29$ in the coupling Hamiltonian \eqref{RabIntHam}.
Furthermore, the external magnetic field parameter is chosen as $C_0=0.35$.
In the classical sector, the initial position and momentum mean values (Input 7) are chosen as $(\mu_q,\mu_p)=(0,4)$.
Note that the chosen parameters result in the initial kinetic energy $E_{\operatorname{kin}}(0)=8$.
Moreover, in the quantum sector, the initial density matrix represents the excited state, which corresponds to the choice $\hat\rho_0=v_0v_0^\dagger$ for $v_0=(1,1)^T/\sqrt{2}$.

\paragraph{Classical sector.}
\begin{figure}[h]
	\begin{center}
		\addtolength{\leftskip}{-20mm}
		\addtolength{\rightskip}{-20mm}
		\includegraphics[width=172mm]{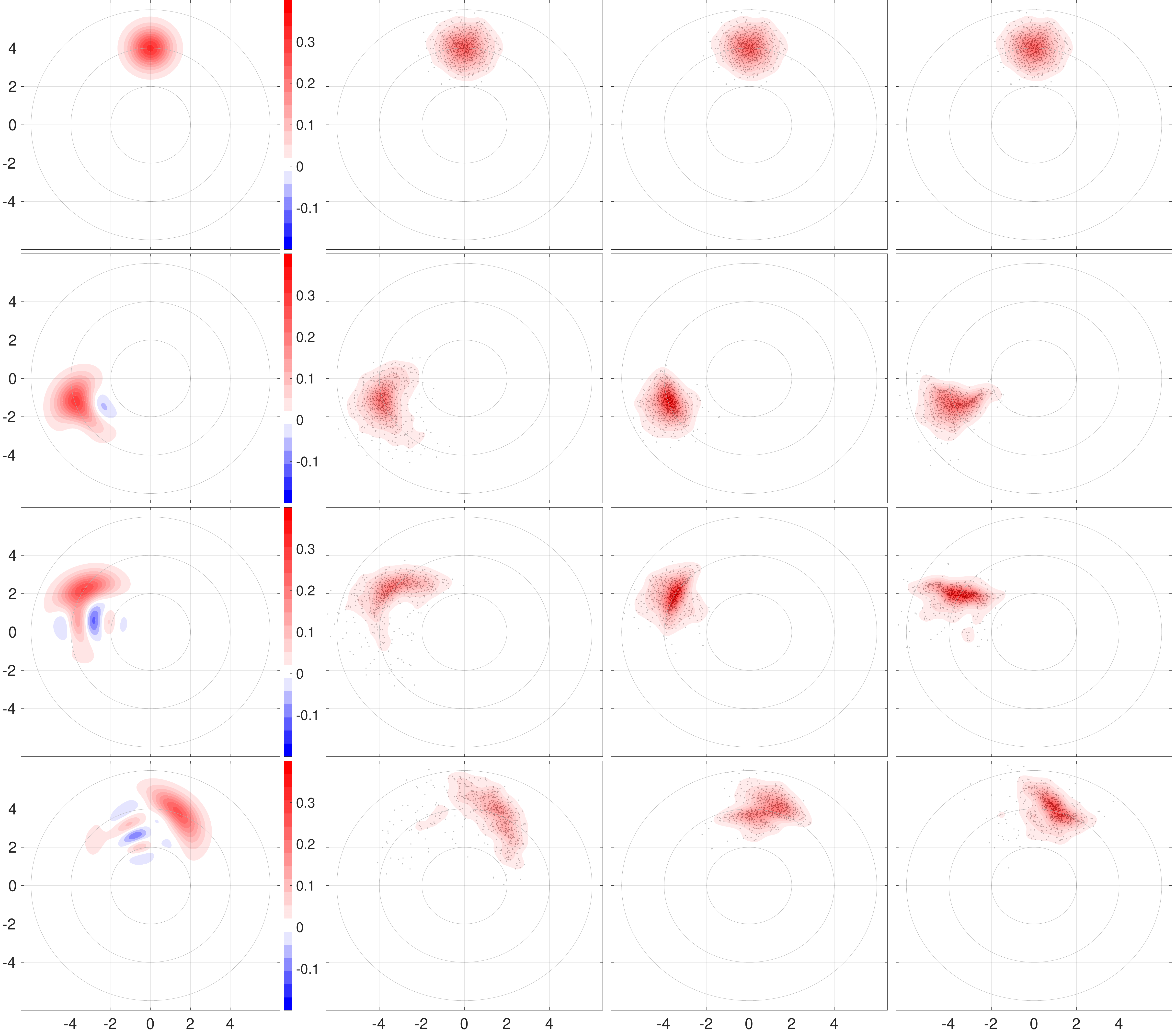}\vspace{-.25cm}
		\caption{\footnotesize
			Time evolution in the classical sector for the Rabi model in the ultrastrong coupling regime.
			Each column (1st: SOFT, 2nd: koopmons, 3rd: Ehrenfest, 4th: bohmions) shows a snapshot of the distribution in phase space.
			Momentum on the vertical axis, position on the horizontal axis (both in atomic units).
			The rows correspond to times $t=0$, $10.5$, $17.5$, and $25$ (top to bottom, all in atomic units).
			Also, $N=500$ and $\alpha=0.5$.\vspace{-5.5mm}
		}\label{US_c}
	\end{center}
\end{figure}
Figure~\ref{US_c} shows the snapshots up to the final time $t=25$ of the simulation.
Equipped with a positive initial momentum, the wavepacket of the quantum solution undergoes a clockwise circular motion, completing four oscillations.
Over time, the shape of the Wigner distribution deviates from its original Gaussian shape.
While both Ehrenfest and bohmion simulations fall short in reproducing the correct shape of the time-evolved wavepacket, the koopmons show a remarkable capability to capture even fine details of the positive regions of the Wigner function.
This includes nuances like the splitting of the tail in the snapshot at the final time.
\begin{figure}[h]
	\centering
	\includegraphics[width=\textwidth]{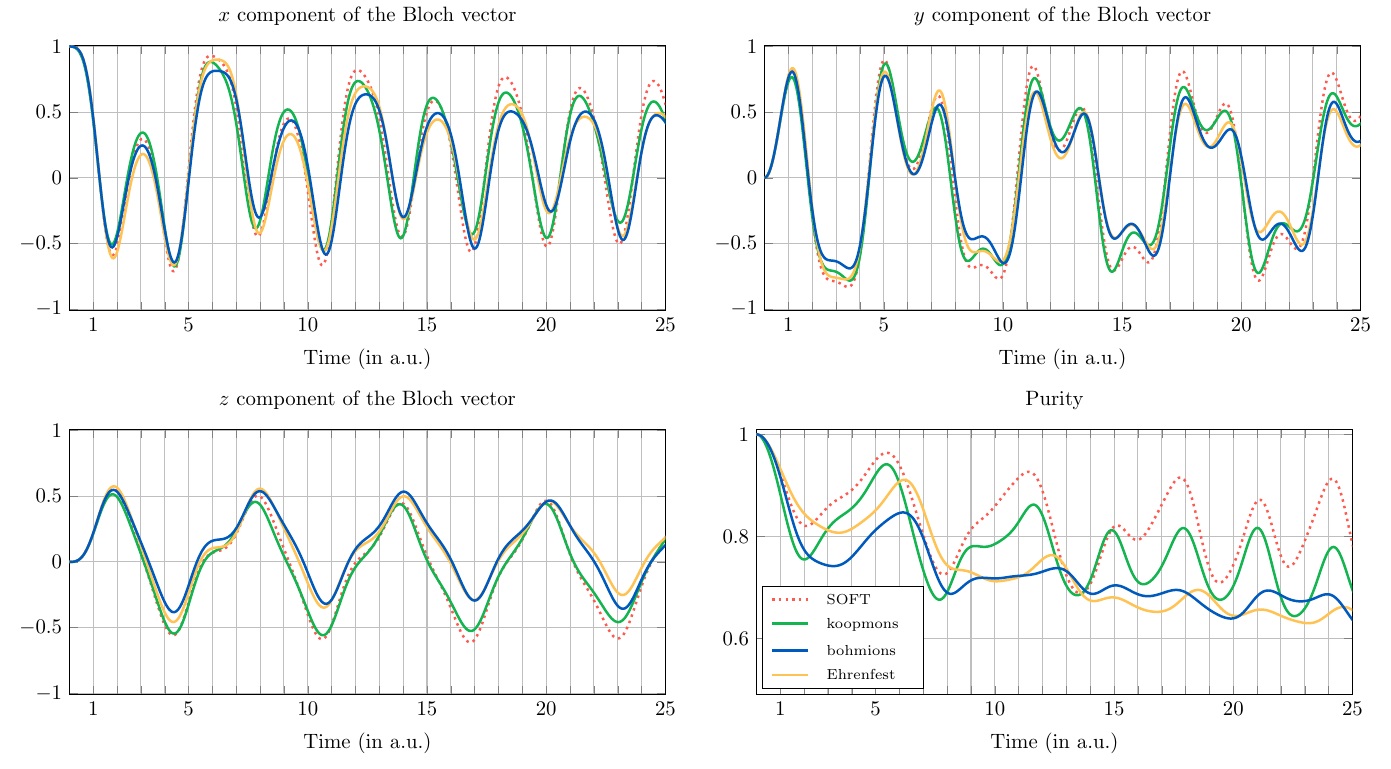}\vspace{-.5cm}
	\caption{\footnotesize
		Components of the bloch vector and quantum purity (bottom left) for the Rabi model in the ultrastrong coupling regime.
	}\label{US_q}
\end{figure}

\paragraph{Quantum sector.}
The results in the quantum sector can be found in Figure~\ref{US_q}.
The purity plot (bottom right) clearly indicates that the oscillations of the quantum solution are most accurately captured by the koopmon method.
In particular, a notable correspondence is observed between the oscillation frequency of the quantum solution and that of the koopmon simulation.
In contrast, both the Ehrenfest and bohmion simulations suffer from overdecoherence at larger times.
Finally, we would like to point out that in the purity plot, differences between the koopmons and the quantum solution appear to grow at larger times.
As previously described for the Tully models, these discrepancies seem to be directly linked to the extent of negativities in the Wigner distribution.
Indeed, while negativities were already observed for the Tully models, it is now evident that the negative regions expand over time.
This means that ultrastrong Rabi dynamics can hardly be modeled by MQC methods since the oscillator acquires purely quantum features as time goes by.

In addition to purity, we also plotted the three components of the Bloch vector $\boldsymbol{b}(t)\in\mathbb{R}^3$.
These components can be expressed in terms of the density matrix, as $b_x(t)=2\operatorname{Re}(\hat\rho_{12}(t))$, $b_y(t)=2\operatorname{Im}(\hat\rho_{21}(t))$, and $b_z(t)=\hat\rho_{11}(t)-\hat\rho_{22}(t)$.
A look at the components of the Bloch vector reveals that all methods qualitatively reproduce the pattern produced by the quantum solution.
However, when compared to the koopmon method, both the Ehrenfest and bohmion simulations show difficulties to reproduce the values obtained by the Schr\"odinger solver.

\subsection{Deep strong coupling regime}
As we will see below, this test case represents a severe challenge for all particle schemes considered in this paper, although the koopmons succeed in reproducing much of the phase-space dynamics of the fully quantum solution.
We consider a case of deep strong coupling regime by setting the constant $\gamma=1.85$ in the coupling Hamiltonian \eqref{RabIntHam}.
Moreover, the external magnetic field parameter is chosen as $C_0=0.1$.
A plot of the corresponding PESs can be found at the right panel of Figure~\ref{pes_R}.
In the classical sector, the initial position and momentum mean values (Input 7) are chosen as $(\mu_q,\mu_p)=(0,0)$.
Note that the chosen parameters result in the initial kinetic energy $E_{\operatorname{kin}}(0)=0$, which is right at the intersection between the two PESs; see Figure~\ref{pes_R}.
Moreover, in the quantum sector, the initial density matrix represents the excited state, which corresponds to the choice $\hat\rho_0=v_0v_0^\dagger$ for $v_0=(1,1)^T/\sqrt{2}$.

\paragraph{Classical sector.}
\begin{figure}[h!]
	\begin{center}
		\addtolength{\leftskip}{-20mm}
		\addtolength{\rightskip}{-20mm}
		\includegraphics[width=172mm]{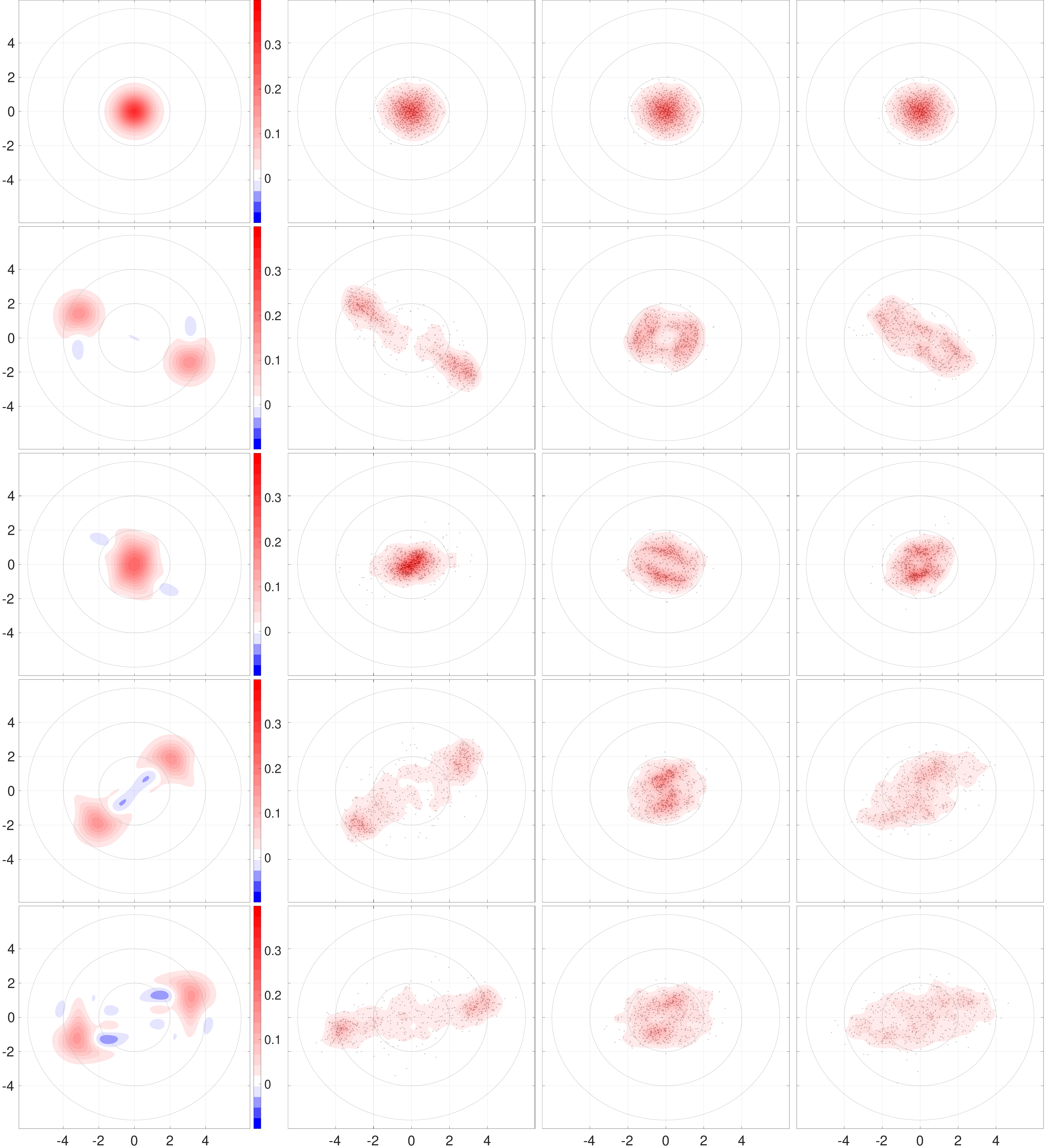}\vspace{-.25cm}
		\caption{\footnotesize
			Time evolution in the classical sector for the Rabi model in the deep strong regime.
			Each column (1st: SOFT, 2nd: koopmons, 3rd: Ehrenfest, 4th: bohmions) shows a snapshot of the distribution in phase space.
			Momentum on the vertical axis, position on the horizontal axis (both in atomic units).
			The rows correspond to times $t=0$, $4$, $6$, $8$ and $15$ (top to bottom, all in atomic units).
			Also, $N=500$ and $\alpha=0.5$.\vspace{-5.5mm}
		}\label{DS_c}
	\end{center}
\end{figure}
Figure~\ref{DS_c} shows the snapshots up to the final time $t=15$ of the simulation.
Starting from the origin in phase space, the Wigner distribution undergoes a splitting into two parts.
By the time $t\approx 3$, the separation between these parts reaches its maximum distance, followed by a reconnection at $t\approx 6$.
Hence, over the entire simulation time, two complete oscillations are captured.

While the Ehrenfest simulation reveals small formations of clusters and rotations in the particle clouds, it fails to reproduce the splitting.
Better results are obtained using the bohmion method.
As we can see, particles get squeezed in a direction that corresponds to the formation of blobs in the quantum solution.
Once more, the koopmons yield the best result.
Especially at early times, we can clearly see the formation of two clusters that carry most of its mass in the right place.
At some point, these clusters reconnect and then split again.
As for the ultrastrong coupling, we can see that the Wigner distribution develops negative values.
However, in the deep strong coupling regime, the negative values not only appear earlier in time, but also rapidly take up larger phase-space regions.

\begin{figure}
	\centering
	\includegraphics[width=\textwidth]{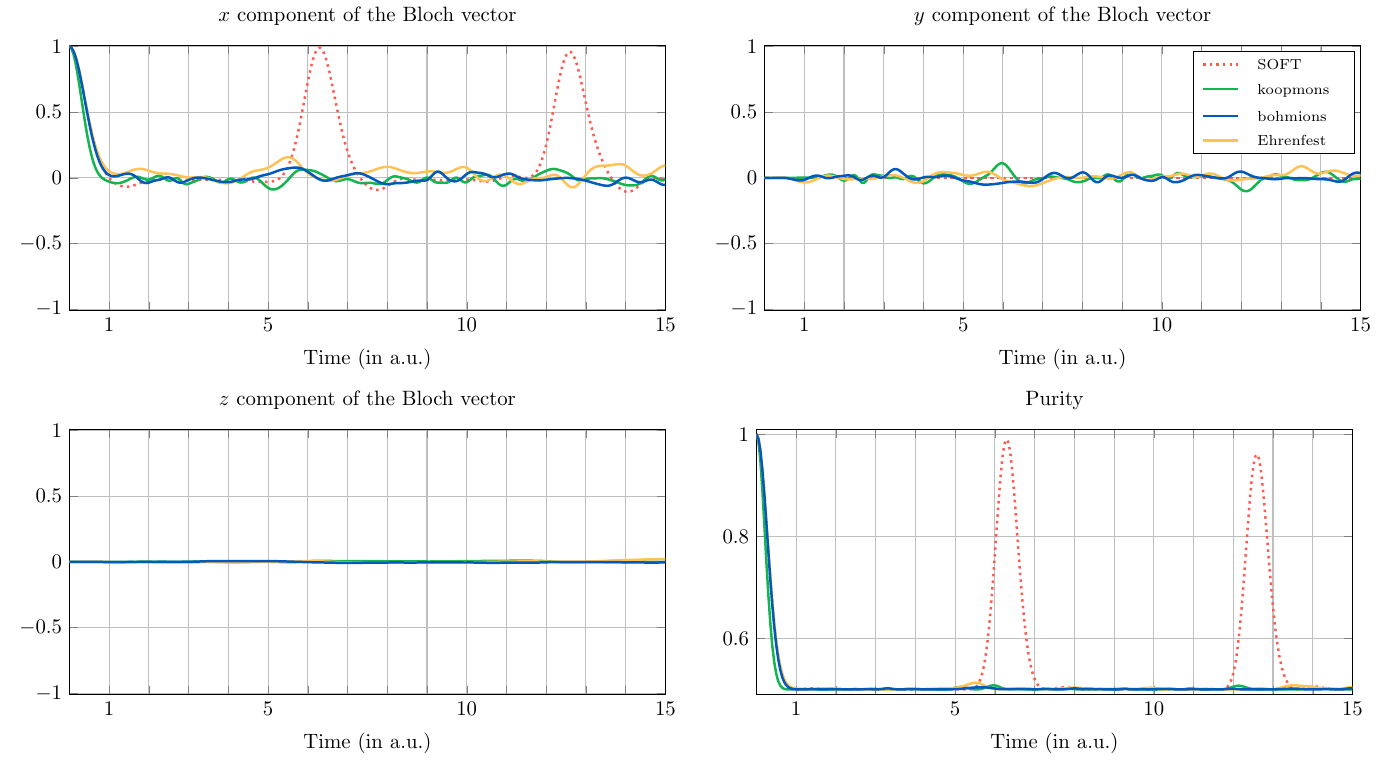}\vspace{-.5cm}
	\caption{\footnotesize
		Components of the bloch vector and quantum purity (bottom left) for the Rabi model in the deep strong coupling regime.
	}\label{DS_q}
\end{figure}

\paragraph{Quantum sector.}
The results in the quantum sector can be found in Figure~\ref{DS_q}.
The quantum solution shows dynamics solely in the $x$-component of the Bloch vector.
None of the considered particle methods can capture this dynamics.
While all particle methods correctly capture the absence of dynamics in the $z$-component, small oscillations are observed in the $y$-component for all methods.
By looking at the purity plots, we can see that all particle methods reproduce the initial decay but fail to capture the revivals around times $t\approx 6$ and $t\approx 12.5$.
The presence of negativities in the Wigner distribution in classical space provides a hint as to why these revivals are elusive:
as mentioned earlier, regions with negative values not only manifest very early in time, but also expand in size.
This means that, similarly to ultrastrong Rabi dynamics, the deep strong regime can hardly be modeled by MQC methods since the oscillator rapidly acquires purely quantum features that escape the MQC treatment.
In particular, we would like to note that further tests have demonstrated that decreasing the regularization parameter $\alpha$ or increasing the number of particles does not help in capturing the revival of purity.

\subsection{Dynamics in configuration space}
So far, we have shown that the koopmons do better than the bohmions in reproducing the fully quantum solution for the present values of the regularization parameter $\alpha$.
Part of this comparison was done by considering the dynamics in phase space.
However, one may reasonably argue that this is not a fair comparison because the bohmion scheme arises in the context of quantum hydrodynamics and thus only densities in configuration space should be considered, rather than dynamics on phase space.
In order to address this point, in this section we present a comparison study that is entirely based on the configuration space and focuses especially on Rabi dynamics in the ultrastrong and deep strong regimes.
In particular, we present \emph{waterfall plots} of the density profile in coordinate space at different times.
\begin{figure}[t]
	\centering
	\includegraphics[width=17.25cm]{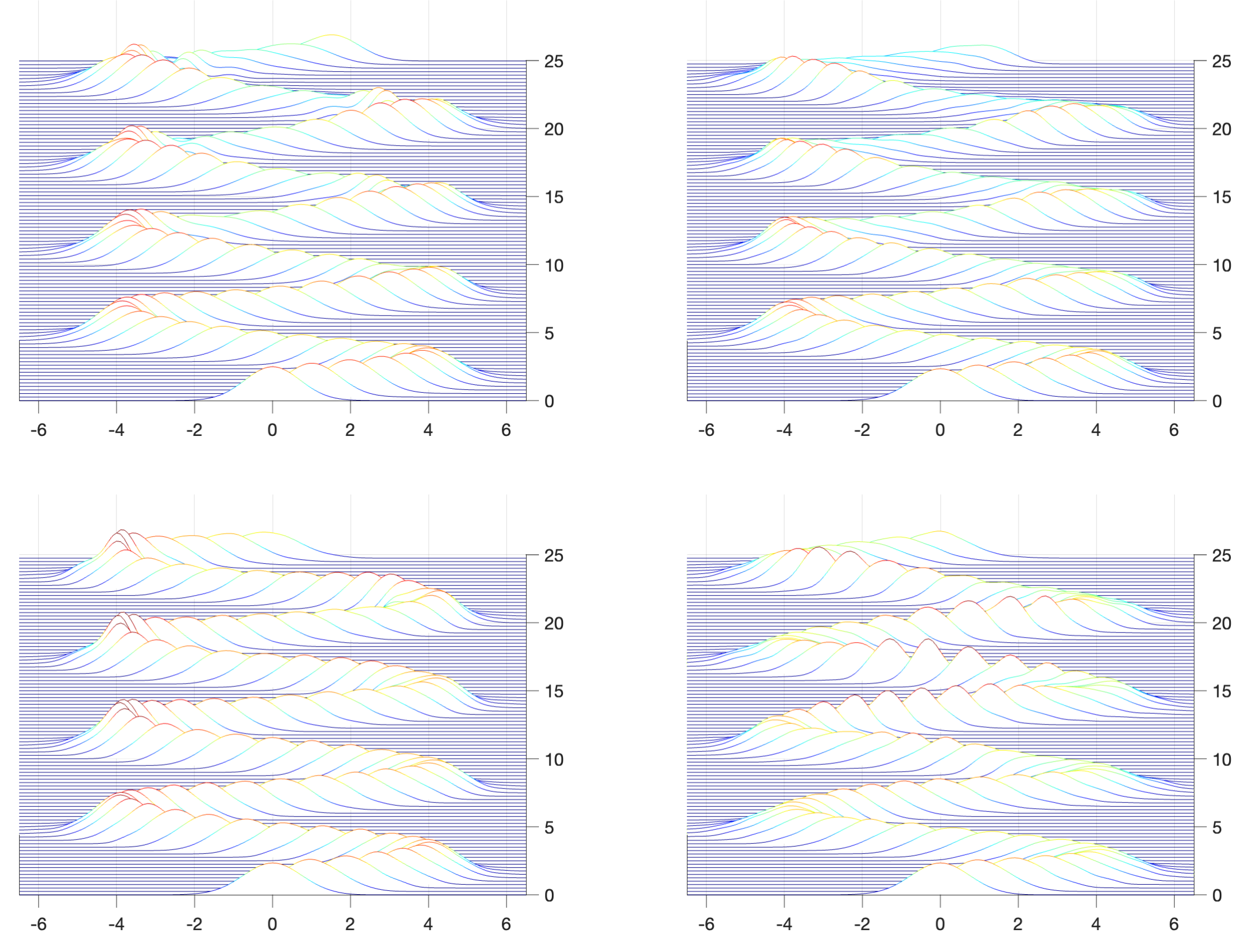}\vspace{-.5cm}
	\caption{\footnotesize
		Waterfall plots of the density in configuration space for the Rabi model in the ultrastrong coupling regime.
		Top: SOFT (left), koopmons (right).
		Bottom: Ehrenfest (left), bohmions (right).
		Time on the vertical axis, position on the horizontal axis (both in atomic units).
		In addition, $N=500$ and $\alpha=0.5$.
	}\label{US_w}
\end{figure}
For the SOFT simulation we plot $|\Psi(r,t)|^2$, while for the particle methods, we plot the density $\tilde D(r,t)=\sum_a w_a\tilde K^{(\Delta)}(r-q_a(t))$ for different times $t$.
Here, the visualization kernel $K^{(\Delta)}$ was introduced in Section~\ref{sub:Quantum reference solver and visualization}.
In the case of ultrastrong coupling, these plots can be found in Figure~\ref{US_w}.
Across all panels it is evident that the characteristic oscillations are captured by all particle methods.
However, another distinctive pattern can be observed in the quantum solution (top left):
As the wavepacket approaches its maximum and minimum position values, it gets more localized, whereas the spread increases in between.
This behaviour is particularly visible after the first three oscillations and gets more pronounced over time.
For example, observing the maximal heights of the wavepacket in the last two oscillations reveals a rapid transition in the spread.
While the bohmions (bottom right) clearly have difficulties to reproduce this pattern, the Ehrenfest simulation (bottom left) yields a qualitatively better result.
The most favourable outcomes, however, are achieved by the koopmons (top right).
\begin{figure}[h]
	\centering
	\includegraphics[width=17.25cm]{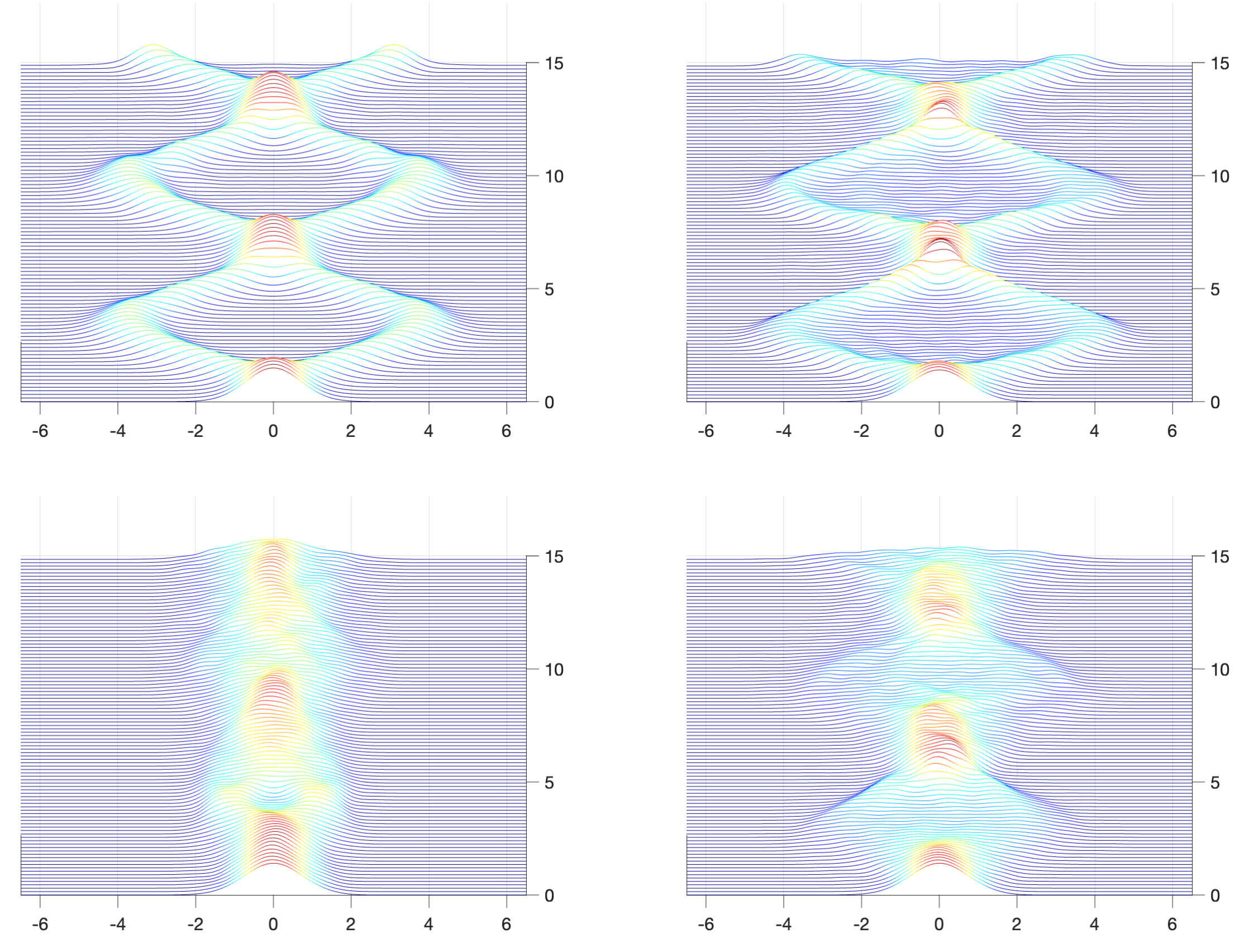}\vspace{-.5cm}
	\caption{\footnotesize
		Waterfall plots of the density in configuration space for the Rabi model in the deep strong coupling regime.
		Top: SOFT (left), koopmons (right).
		Bottom: Ehrenfest (left), bohmions (right).
		Time on the vertical axis, position on the horizontal axis (both in atomic units).
		In addition, $N=500$ and $\alpha=0.5$.
	}\label{DS_w}
\end{figure}

In the case of deep strong coupling, the waterfall plots can be found in Figure~\ref{DS_w}.
We can see that the branching of the wavepacket is accurately captured by the koopmon method.
Although the reconnection times of the branches in the bohmion method match those of the quantum solution, the branches themselves have a smaller maximum span.
The weakest results are obtained with the Ehrenfest method, where almost no branching is visible.

\section{Koopmons and bohmions in higher dimensions\label{sec:Extensions to higher dimensions}}
As discussed in Section~\ref{sec:Regularization and the particle scheme}, the trajectory-coupling in the koopmon method relies on the computation of the multi-dimensional integrals $\widehat{\cal I}_{ab}$, which model the quantum-classical correlation effects.
In fact, when dealing with a classical subsystem of dimension $d\ge1$ in configuration space, numerical integration must be executed over $\mathbb{R}^{2d}$.
Consequently, conventional grid-based approaches become impractical as the number of grid points grows exponentially with the dimension, leading to unacceptably slow computations.
While advanced integration techniques like sparse grids \cite{Gerstner98} or low-rank approximations \cite{Oseledets11} can help to overcome the curse of dimensionality to some extent, they still fall short for high-dimensional model systems, even when high-performance computing facilities are used.

To make the koopmon method applicable to relevant situations in chemistry and physics, avoiding high-dimensional integrals in the model equations becomes a crucial requirement.
In this section, we demonstrate that, under rather general conditions on the form of the Hamiltonian $\widehat{H}$, the koopmon method can be implemented without the need to compute integrals of dimensions higher than two.
In this section, we present our results for two degrees of freedom.
The extension to higher dimensional cases is straightforward.

Consider two degrees of freedom with phase-space coordinates $\bz_1=(q_1,p_1)$ and $\bz_2=(q_2,p_2)$.
For a given multi-index $\kappa\in\{(a,b)\,:\,a,b=1,\dots,N\}$, we adopt the procedure from Section~\ref{sec:Regularization and the particle scheme} and assume the singular solution ansatz in two dimensions, that is,
\begin{align*}
	\widehat{\cal P}(\bz_1,\bz_2,t)
	=\sum_{\kappa} w_\kappa\hat\varrho_\kappa(t)\delta(\bz_1-\bzeta_1^{(\kappa)}(t))\delta(\bz_2-\bzeta_2^{(\kappa)}(t)).	
\end{align*}
In addition, we restrict this ansatz to consider the case $w_\kappa=w_{(a,b)}=w_a w_b$, along with
\beq\label{MIansatz}
	\bzeta_1^{(a,b)}
	=\bzeta_1^{(a,b')}\qquad\forall a,b,b'
	\qquad\text{and}\qquad
	\bzeta_2^{(a,b)}
	=\bzeta_2^{(a',b)}\qquad\forall b,a,a'.
\eeq
Furthermore, we pick an even kernel $K\in C^2(\mathbb{R}^4)$ such that $K(\bz-\bz')=K_1(\bz_1-\bz_1')K_2(\bz_2-\bz_2')$, so that
\begin{align*}
	\bar{\cal P}(\bz_1,\bz_2,t)
	&=\sum_{a,b}w_a w_b\hat\varrho_{ab}(t)K_1(\bz_1-\bzeta_1^{(a)}(t))K_2(\bz_2-\bzeta_2^{b}(t))\\
	&=:\sum_{a,b}w_aw_b\hat\varrho_{ab}(t)K_1^{(a)}(\bz_1,t)K_2^{(b)}(\bz_2,t),
\end{align*}
as well as
$
	\bar{f}
	=\operatorname{Tr}\bar{\cal P}
	=\sum_{a,b} w_a w_b K_1^{(a)}K_2^{(b)}
	=\big(\sum_{a} w_aK_1^{(a)}\big)\big(\sum_{b}w_bK_2^{(b)}\big)
$.
We compute
\begin{align*}
	&\frac12\sum_{\kappa,\kappa'}w_\kappa w_{\kappa'}\bigg\langle i\hbar\big[\hat{\varrho}_{\kappa'},\hat{\varrho}_\kappa\big],\iint\frac{K_\kappa\{K_{\kappa'},\widehat{H}\}}{\sum_\gamma w_\gamma K_\gamma}\,\de^2z_1\de^2z_2\bigg\rangle\\
	&\qquad\qquad
	=\frac12\sum_{a,b,a',b'}w_aw_bw_{a'}w_{b'}\bigg\langle i\hbar\big[\hat{\varrho}_{a'b'},\hat{\varrho}_{ab}\big],\,\widehat{\!{I}}_{\!aba'b'},\bigg\rangle\\
	&\qquad\qquad=\frac14\sum_{a,b,a',b'}w_aw_bw_{a'}w_{b'}\bigg\langle i\hbar\big[\hat{\varrho}_{a'b'},\hat{\varrho}_{ab}\big],\,\widehat{\!{I}}_{\!aba'b'}-\widehat{\!{I}}_{\!a'b'ab}\bigg\rangle,
\end{align*}
where we have denoted
\begin{align*}
	\widehat{\!{I}}_{\!aba'b'}
	:=\iint\bigg(\frac{K_1^{(a)}K_2^{(b)}K_2^{(b')}\{K_1^{(a')},\widehat{H}\}_1}{\sum_{c,d} w_cw_d K_1^{(c)}K_2^{(d)}}+\frac{K_1^{(a)}K_2^{(b)}K_1^{(a')}\{K_2^{(b')},\widehat{H}\}_2}{\sum_{c,d} w_cw_d K_1^{(c)}K_2^{(d)}}\bigg)\,\de^2z_1\de^2z_2.
\end{align*}
Hence, by re-organizing the terms of the integrand and applying Fubini's theorem, we get
\begin{align*}
	\widehat{\!{I}}_{\!aba'b'}
	&=\iint\bigg(\frac{K_2^{(b)}K_2^{(b')}}{\sum_{d} w_d K_2^{(d)}}\frac{K_1^{(a)}\{K_1^{(a')},\widehat{H}\}_1}{\sum_{c} w_c K_1^{(c)}}+\frac{K_1^{(a)}K_1^{(a')}}{\sum_{c}w_c K_1^{(c)}}\frac{K_2^{(b)}\{K_2^{(b')},\widehat{H}\}_2}{\sum_{d}w_dK_2^{(d)}}\bigg)\,\de^2z_1\de^2z_2\\
	&=\int\!\frac{K_1^{(a)}}{\sum_{c}w_cK_1^{(c)}}\bigg\{K_1^{(a')},\int\frac{K_2^{(b)}K_2^{(b')}\widehat{H}}{\sum_{d} w_d K_2^{(d)}}\,\de^2 z_2\bigg\}_1\!\de^2z_1\\
	&\qquad\qquad\qquad+\int\!\frac{K_2^{(b)}}{\sum_{d} w_dK_2^{(d)}}\left\{K_2^{(b')},\int\frac{K_1^{(a)}K_1^{(a')}\widehat{H}}{\sum_{c} w_c K_1^{(c)}}\de^2z_1\right\}_2\!\de^2z_2.
\end{align*}
Without further assumptions, no simplifications of the last two integrals can be made.
Consequently, we are confronted with the integration in four dimensions.
However, as we will see below, a rather nonrestrictive assumption on the form of the Hamiltonian allows a factorization into two-dimensional integrals.
In particular, we consider the following class of Hamiltonians:
\begin{align}\label{classofH}
	\widehat{H}(\bz_1,\bz_2)
	=H_c(\bz_1,\bz_2)\boldsymbol{1}+\widehat{H}_Q+h_1(\bz_1)\widehat{H}_2(\bz_2)+h_2(\bz_2)\widehat{H}_1(\bz_1),
\end{align}
which, to the best of our knowledge, covers most common models in computational chemistry.
A short calculation shows that
\begin{align*}
	\widehat{\!{I}}_{\!aba'b'}
	=\widehat{\cal I}_{aba'b'}+C_{aba'b'}\boldsymbol{1},
\end{align*}
where $C_{aba'b'}$ is a real-valued constant, and
\begin{align*}
	\widehat{\cal I}_{aba'b'}
	=\widehat{\cal I}_1^{(aa')}{\cal J}_2^{(bb')}+{\cal I}_1^{(aa')}\widehat{\cal J}_2^{(bb')}+{\cal I}_2^{(bb')}\widehat{\cal J}_1^{(aa')}+\widehat{\cal I}_2^{(bb')}{\cal J}_1^{(aa')}	
\end{align*}
along with the definitions
\begin{align*}
	{\cal I}_\ell^{(ss')}
	:=\int\frac{K_\ell^{(s)}\{K_\ell^{(s')},h_\ell\}_\ell}{\sum_{c}w_c K_\ell^{(c)}}\,\de^2z_\ell,\qquad
	{\cal J}_\ell^{(ss')}
	:=\int\frac{K_\ell^{(s)}K_\ell^{(s')}h_\ell}{\sum_{d} w_d K_\ell^{(d)}}\,\de^2 z_\ell,	
\end{align*}
\begin{align*}
	\widehat{\cal I}_\ell^{(ss')}
	:=\int\frac{K_\ell^{(s)}\{K_\ell^{(s')},\widehat{H}_\ell\}_\ell}{\sum_{c}w_c K_\ell^{(c)}}\,\de^2z_\ell,\qquad
	\widehat{\cal J}_\ell^{(ss')}
	:=\int\frac{K_\ell^{(s)}K_\ell^{(s')}\widehat{H}_\ell}{\sum_{d}w_d K_\ell^{(d)}}\,\de^2 z_\ell,	
\end{align*}
for $\ell=1,2$ and $s=a,b$.
Finally, using that $\langle C_{aba'b'}\boldsymbol{1},i\hbar[\hat{\varrho}_{ab},\hat{\rho}_{a'b'}]\rangle=0$, we conclude that
\begin{align*}
	\Big\langle i\hbar\big[\hat{\varrho}_{ab},\hat{\varrho}_{a'b'}\big],\,\widehat{\!{I}}_{\!aba'b'}-\widehat{\!{I}}_{\!a'b'ab}\Big\rangle
	=\Big\langle i\hbar\big[\hat{\varrho}_{ab},\hat{\varrho}_{a'b'}\big],\,\widehat{\cal I}_{aba'b'}-\widehat{\cal I}_{a'b'ab}\Big\rangle.
\end{align*}
Given that the last term only involves two-dimensional integrals, the above calculations show that the koopmon scheme can be implemented for two-dimensional systems by computing only sums of products of integrals, each of dimension two.
Since the presented factorization technique can be adapted for systems of dimension $d>2$ in configuration space, the koopmon method can be implemented without the need to compute integrals of dimensions higher than two.
Given the generality of the underlying form we require for $\widehat{H}$, we expect that the koopmon method can be applied to several relevant models in the future.
\begin{remark}
	The class of Hamiltonians in \eqref{classofH} can be extended to cases where the factors $h_j(\bz_j),\,j=1,2$, also carry quantum degrees of freedom.
	More precisely, given that the full Hamiltonian $\widehat{H}(\bz_1,\bz_2)$ must be Hermitian, one could, for instance, replace the term $h_1(\bz_1)\widehat{H}_2(\bz_2)$ with the anticommutator $[\hat h_1(\bz_1),\widehat{H}_2(\bz_2)]_+=\hat h_1(\bz_1)\widehat{H}_2(\bz_2)+\widehat{H}_2(\bz_2)\hat h_1(\bz_1)$.
	A straightforward calculation shows that, also for these more general Hamiltonians, the multidimensional integrals $\widehat{\cal I}_{aba'b'}$ can still be decomposed into sums of products of integrals, each of dimension two.
\end{remark}
We note that a similar factorization is applicable to the bohmion method.
Importantly, due to the independence of the coupling integral on the Hamiltonian, the extension to higher dimensions does not require any restrictions and is more straightforward.
Upon restricting again to the case of two degrees of freedom for simplicity, we resort to a multi-index notation and mimic the relations \eqref{MIansatz} so that ${q_1^{(a,b)}=q_1^{(a,b')}}$ and ${q_2^{(a,b)}=q_2^{(a',b)}}$, along with $w_\kappa=w_{(a,b)}=w_a w_b$.
Then, the two-dimensional extension $\widehat{\!\mathscr{I}}_{\!\kappa\kappa'}$ of the bohmion integral \eqref{bohmintegs} can be rewritten in such a way that
\begin{align*}
	\sum_{\kappa\kappa'}w_\kappa w_{\kappa'}\left(2\langle\hat\varrho_{\kappa},\hat\varrho_{\kappa'}\rangle-1\right)\widehat{\!\mathscr{I}}_{\!\kappa\kappa'}=\sum_{aba'b'}w_aw_bw_{a'}w_{b'}\left(2\langle\hat\varrho_{ab},\hat\varrho_{a'b'}\rangle-1\right)\big({\cal I}_{1}^{(aa')}{\cal J}_{2}^{(bb')}+{\cal I}_{2}^{(bb')}{\cal J}_{1}^{(aa')}\big),	
\end{align*}
where we have introduced the one-dimensional integrals
\begin{align*}
	{\cal I}_{\ell}^{(ss')}
	:=\int\frac{\partial_\ell {\sf K}_{\ell}^{(s)}{{\partial_\ell}}{\sf K}_{\ell}^{(s')}}{\sum_{c}w_cK_{\ell}^{(c)}}\de{{r_\ell}}
	\qquad\text{and}\qquad
	{\cal J}_\ell^{(ss')}
	:=\int\frac{{\sf K}_{\ell}^{(s)}{\sf K}_{\ell}^{(s')}}{\sum_{d}w_d{\sf K}_{\ell}^{(d)}}\de r_\ell
\end{align*}
for $\ell=1,2$ and $s=a,b$.
Here, ${\sf K}_{\ell}^{(s)}={\sf K}_\ell(x_\ell-q_\ell^{s})$.
Therefore, the bohmion method only requires integrals of dimension not higher than one, regardless of the Hamiltonian.
In addition, the bohmion approach requires the evaluation of a smaller number of integrals at each time-step.
On the one hand, this may suggest that the bohmion approach may be more advantageous in higher dimensional problems.
On the other hand, we should not forget that the general grid-size restrictions arising from the necessary low values of the parameter $\alpha$ make the bohmion integrals more complex and computationally costly.
The question whether koopmons or bohmions are better suited for higher dimensional problems is left open for further research.

\section{Conclusions and Outlook}\label{sec:Conclusions and Outlook}
In this paper, we presented the formulation of a new MQC numerical algorithm, the \emph{koopmon method}, based on particle trajectories that arise from a statistical sampling of the Lagrangian paths in phase space.
Hinging on the variational formulation of its underlying continuum model, the method retains a series of consistency issues commonly posing problems in conventional approaches.
In addition, the Hamiltonian structure allows to retain basic conservation laws and also leads to capturing quantum decoherence better than other methods.
We emphasize that decoherence is a crucial property which remains elusive in many conventional methods, such as surface hopping and Ehrenfest dynamics.
The method involves the calculation of a phase-space integral at each time step, requiring an integration grid whose size is chosen depending on a modeling parameter.
As the latter tends to zero, the numerical scheme approaches the exact solution of the underlying MQC continuum model.

In order to illustrate the koopmon method, its numerical scheme was implemented on a laptop machine and applied to different test cases in nonadiabatic dynamics.
In addition, the method was compared to the fully quantum Schr\"odinger solution as well as both the MQC Ehrenfest scheme and the bohmion method, an analogous particle based method for fully quantum hydrodynamics.
In the considered parameter range, the koopmon method succeeds in capturing accuracy levels that are not achieved by either Ehrenfest or bohmion dynamics.
In the case of the Tully models, the koopmon method was most successful for Tully~I and III, while a loss in accuracy is observed for Tully~II at intermediate times.
We argue that this loss of accuracy is due to the emergence of negative regions of the Wigner distribution in the fully quantum dynamics.
These negative Wigner regions signal the occurrence of fully quantum features that challenge the MQC approximation.
Nevertheless, unlike the Ehrenfest scheme, the koopmons recover satisfactory accuracy levels after the negative regions fade away.
In particular, we wish to emphasize that the koopmons are particularly successful for Tully~III, which is a particularly challenging problem for common methods.
In particular, for the given set of parameters, the koopmons are the only method that can capture the purity revival over time.

Two more test cases were studied for the Rabi problem: the ultrastrong and deep strong coupling regimes.
In both these cases, the negative regions in the fully quantum Wigner dynamics form and persist indefinitely over time.
Thus, these coupling regimes represent particularly challenging test cases for MQC methods.
Yet, as we have shown, the koopmons succeed to capture the phase-space dynamics of the positive Wigner regions with great accuracy.
In the ultrastrong coupling regime, the koopmons succeed in capturing also the purity levels better than the other particle methods.
Instead, the deep strong coupling represents an equally challenging test case for all particle methods when it comes to reproducing purity dynamics.
We argue that this is due to the fact that the negative Wigner regions develop very quickly and expand significantly over time, thereby challenging the MQC approximation.

As mentioned at the beginning of Section~\ref{sec:Results for the Tully models}, our numerical studies showed that the choice of $\alpha=0.5$ for the regularization parameter, together with $N=500$ for the number of particles, allows for good comparison with the fully quantum results in all test cases.
However, to achieve satisfactory levels of accuracy in both the quantum and classical sectors for the Tully models, adjustments were made to both $\alpha$ and $N$.
This led to the final choice of $\alpha=0.325$ and $N=1000$ for the Tully models.
For a given problem under consideration, a universal prescription for the choice of model parameters remains an open problem and requires a further study of regularization models in general.
Nevertheless, all our current studies tend to confirm that $\alpha=0.5$ and $N=500$ also provide good results for other models.
Although we anticipate that a time-adaptive tuning of the parameters $\alpha$ and $N$ is an interesting topic for numerical studies in the future, this has not been considered in the present work as this step comes at the expense of exact energy conservation and momentum balance.

Our comparison between koopmons and bohmions motivated us to also look at their features in configuration space.
This is motivated by the fact that the bohmions possess an underlying continuum description entirely based on quantum hydrodynamic equations in configuration space, rather than phase space.
An analysis of the density dynamics in configuration space for koopmons, bohmions, and Ehrenfest trajectories essentially confirmed the results already appeared in phase space, with the koopmons achieving higher accuracy in all cases.
Once more, we emphasize that the bohmion method generally allows for much more accurate results than those presented here \cite{HoRaTr21}.
However, increasing the accuracy requires both a considerably finer integration grid and a higher number of particles, thereby increasing the numerical costs.

The final discussion focused on the higher-dimensional extensions of both the koopmon and the bohmion methods.
As we showed, a further closure at the multi-index level allowed for reducing the higher-dimensional integral occurring in both approaches to a quadratic expression involving only two-dimensional and one-dimensional integrals, respectively for the koopmons and the bohmions.
While this lower-dimensional reduction of the integrals appear for the bohmions in all cases, for the koopmons one needs to assume a certain form of the Hamiltonian, which however comprises most current models in computational chemistry.
Thus, we have good reasons to believe that the treatment of higher-dimensional problems will not be computationally prohibitive as it might be without using the convenient multi-index closure.

The higher-dimensional extension will be the subject of future activities, with special focus on the effects of conical intersections between energy surfaces.
The latter continue to represent an active area of research and we would like to understand how the koopmons can help in this direction.
Within the context of nonadiabatic dynamics, we are currently considering also extended Tully models such as the Double Arch model proposed in \cite{Subotnik11}.
In addition, we will also add momentum coupling to approach problems in spintronics and solid-state physics, where the Rashba and Dresselhaus effects play a prominent role in spin-orbit coupling.
Finally, the present comparison of bohmions and koopmons also opens a new interesting direction: how can one couple the bohmions and the koopmons in an MQC model that retains the full infinite-dimensional character of the dynamics in the quantum sector?
These and related questions will be explored in future work.

Regarding the numerical implementation of the koopmon method, we identify potential improvements that can help to reduce the numerical costs or increase the accuracy.
Embedded RK methods such as the celebrated Dormand--Prince 4(5) (RKDP) or Fehlberg's 4(5) (RKF) method enable efficient implementation and local error control with adaptive step sizes.
Both the RKDP and RKF methods have recently been added to our code.
Furthermore, since the underlying model equations are Hamiltonian, a symplectic time integrator should be advantageous.
We point out that the underlying Hamiltonian function is not generally separable and therefore explicit symplectic integrators are not available or rely on a complicated construction in double phase space \cite{Tao}.
While implicit symplectic RK methods come at the expense of solving a nonlinear system, we have recently implemented the fourth-order symplectic two-stage RK method based on Gaussian quadrature \cite{Hairer}.
In this case, we solve the underlying nonlinear system with fixed-point iteration.

The numerical integration of the backreaction integrals, which are currently computed using the composite trapezoidal rule, could also be improved.
Based on the choice of Gaussian kernel functions, two-dimensional Gauss--Hermite quadrature could be considered for future implementation.
Alternatively, discretization of the phase-space integrals using Monte Carlo or quasi-Monte Carlo quadrature offers a promising alternative which have already been successfully implemented in the context of high-dimensional applications of quantum molecular dynamics \cite{LasserHK}.
Finally, it should be noted that the equations of the koopmon method allow for the parallel computation of the backreaction integrals for each trajectory due to the product structure $K(q,p)=\tilde K(q)\tilde K(p)$ of the kernel functions.
As we have seen for the Tully III model, the koopmons can span a considerable range of phase space at longer simulation times, making the current choice of the spatial grid for the composite trapezoidal rule a bottleneck.
The parallel computation of the phase-space integrals with corresponding local grids holds a significant optimization potential for the current implementation, which we also plan to explore in future work.

Another future direction also involves the possibility of a so-called \emph{on-the-fly} implementation of koopmon dynamics.
This direction is motivated by the fact that in many practical cases the MQC Hamiltonian $\widehat{H}(\bz)$ is not known at all phase-space points, while it is only computed along trajectories.
For example, in our case this may be achieved by a modified regularization process in such a way that only the quantities $\nabla\widehat{\cal P}$ and $f$ are convoluted in the backreaction integrand $\langle i\hbar\{\widehat{\cal P},\widehat{H}\}\rangle={f^{-1}\operatorname{Tr}(\bX_{\widehat{H}}\cdot[i\hbar\widehat{\cal P},\nabla\widehat{\cal P}])}/2$ appearing in \eqref{VarPrin1}.
While this approach treats $\widehat{\cal P}$ and $\nabla\widehat{\cal P}$ on a different footing, it may be worth pursuing to make further progress.
This and other options will be considered elsewhere.

\bigskip
\paragraph{Acknowledgments.}
{\small
We are all greatly indebted to Darryl Holm, who is part of the collaboration working on the bohmion scheme and provided extremely valuable advice during our discussion meetings with him at Imperial College London.
CT also wishes to thank Robert MacKay for his encouragement to work in this direction and for his insistence on the need to apply the new model to problems in physics and chemistry.
It turned out very rewarding advice.
We also thank all the other people who supported our efforts and provided further motivation during this work, especially Victor S. Batista, Michael Berry, Denys Bondar, Irene Burghardt, Basile Curchod, Eberhard Gross, Caroline Lasser, Robert Littlejohn, Giovanni Manfredi, Phil Morrison, Micheline B. Soley, John Tully and Tomasz Tyranowski.
We would also like to thank Delyan Zhelyazov for his recent work on the adaptive-timestep implementation.
Finally, we thank the reviewers for pointing out the importance of nonadiabatic couplings and for drawing our attention to the literature on the Koopman operator.
This work was made possible through the support of Grant 62210 from the John Templeton Foundation.
The opinions expressed in this publication are those of the authors and do not necessarily reflect the views of the John Templeton Foundation.
Financial support by the Leverhulme Research Project Grant RPG-2023-078 is also greatly acknowledged.
FGB was partially supported by a start-up grant from the Nanyang Technological University.
}

\appendix
\section{Role of the nonadiabatic coupling}\label{app:sec:Nonadiabatic coupling}
This Appendix expands on the phenomenology associated to the nonadiabatic processes observed in the Tully models from Section~\ref{sec:Results for the Tully models}.

\begin{figure}[h]
	\centering
	\includegraphics[width=\textwidth]{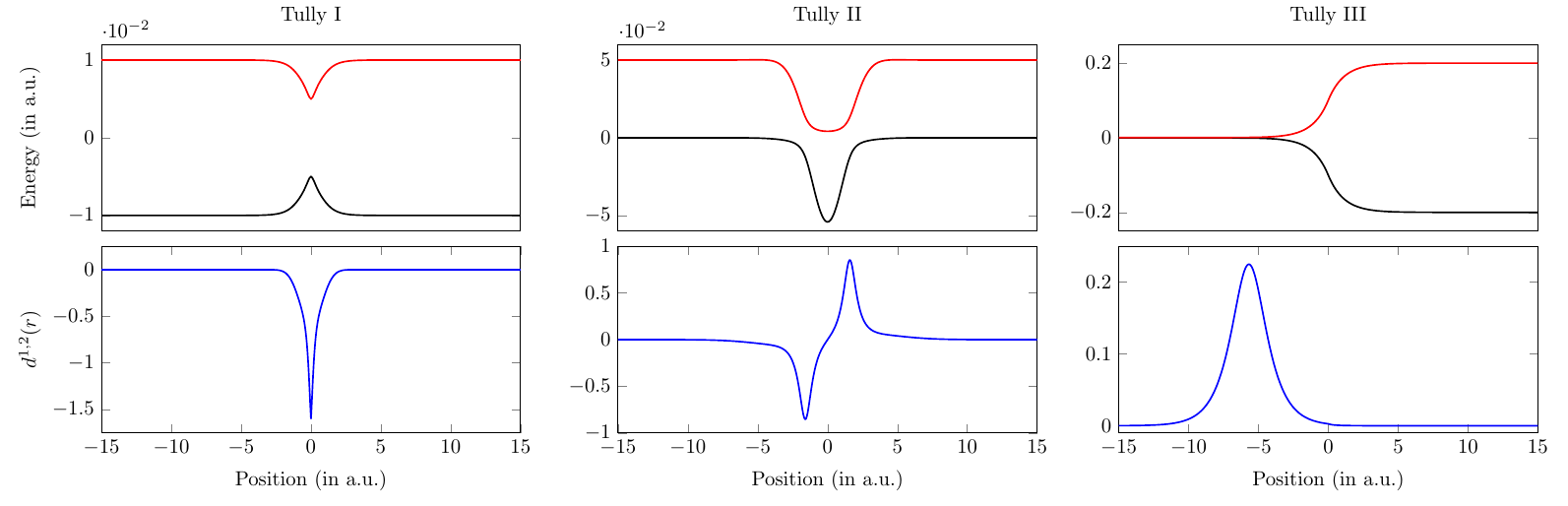}\vspace{-.5cm}
	\caption{\footnotesize
		PESs (top) and nonadiabatic couplings (bottom) for the Tully models.
	}\label{nac_T}
\end{figure}

In the context of the Tully models, the koopmon trajectories model the interaction between nuclear degrees of freedom (classical subsystem), described by the classical trajectories $(q_a,p_a)$, and electronic degrees of freedom (quantum subsystem), described by the density matrices $\hat\varrho_a$.
In the adiabatic regime, the nuclear dynamics take place on the two-level potential energy surfaces (PESs), which are shown in the top row of Figure~\ref{nac_T}.
These surfaces are sometimes referred to as \emph{adiabatic} PESs, reflecting the fact that they are given by the eigenvalues of the electronic Hamiltonian ($\widehat{H}_{\operatorname{el}}(q)={H_0(q){\boldsymbol{1}}+H_1(q)\hat\sigma_x+H_3(q)\hat\sigma_z}$), with the eigenvectors forming the so-called \emph{adiabatic} basis.
In nonadiabatic dynamics, however, one is usually not only interested in the nuclear dynamics on adiabatic PESs, but also in the exchange of quantum amplitudes between the surfaces.
The typical picture involves the Born--Huang expansion
\begin{align*}
	\Psi(r,t)
	=\psi_1(r,t)v_1(r)+\psi_2(r,t)v_2(r)
\end{align*}
of the molecular wave function $\Psi(\cdot,t)\in L^2(\mathbb{R},\mathbb{C}^2)$ in the adiabatic basis $\{v_1(r),v_2(r)\}\subset\mathbb{C}^2$ made available by the spectral problem for the electronic Hamiltonian.
In this context, the so-called \emph{nonadiabatic couplings} (NACs) are often used to identify regions in the PES landscape where nonadiabatic processes are likely to be observed.
In regions where the PESs come close to each other, so-called \emph{nonadiabatic} transitions occur, and the different Tully models were designed to create a platform for testing models to capture these effects.

NACs are expressed in terms of NAC vectors, which are particularly important for computing the standard jump probability for the fewest switches surface hopping algorithm introduced by Tully \cite{Tully90}.
For two-level systems defined on a one-dimensional coordinate space, the NAC vector (in this case a real number) is given by
\begin{align}\label{NACv}
	d^{1,2}(r)
	:=\langle v_1(r)\mid\nabla v_2(r)\rangle
	=\frac{\langle v_1(r)\mid\nabla H_{\operatorname{el}}(r)v_2(r)\rangle}{\lambda_1(r)-\lambda_2(r)},
\end{align}
where $\lambda_j(r),\,j\in\{1,2\}$, denotes the corresponding eigenvalues (and hence PES) of the adiabatic basis vectors $v_j(r)$ at a given position $r\in\mathbb{R}$.
The last term in \eqref{NACv} provides an alternative way to compute the NAC, and the proof of the second equality goes back to \cite{Born28}.

Figure~\ref{nac_T} shows a plot of the NACs (blue) for the three Tully models.
As we can see, the values of $d^{1,2}(r)$ are indeed large in regions where the PESs (top) are close to each other.
Furthermore, formula~\eqref{NACv} reveals that the NAC vectors have a singularity for positions with $\lambda_1(r)=\lambda_2(r)$, and it is known that the adiabatic approximation breaks down near these points.

As we can see in Figures~\ref{T1_c}, \ref{T2_c} and \ref{T3_c}, regions of strong nonadiabatic coupling cause significant deformations of the initial Gaussian wavepacket entering from the left.
Moreover, the presence of negative values in these regions indicates purely quantum effects.
Once these negativities disappear, the population levels reach a constant value, indicating that there is no further exchange of quantum amplitudes.
For the Tully models I and II, we find that in regions with strong nonadiabatic coupling there is a certain loss of accuracy in reproducing the correct population levels as predicted by the fully quantum code.
However, the koopmons reproduce the correct population levels before entering and after leaving the regions of strong nonadiabatic coupling.
For Tully III, the plot of the NAC shows that the nonadiabatic coupling reaches its maximum around $x\approx-6$.
The koopmons are unique in qualitatively capturing the correct population dynamics for all times (see Figure~\ref{T3_pp}).
Both the Ehrenfest and bohmion methods fail to capture the reflection (no negative momentum for the trajectories), while the koopmons reproduce negative momenta as predicted by the fully quantum simulation.
Since only the reflected trajectories revisit the region of strong nonadiabatic coupling, more transitions between the PESs take place between times $t = 2000$ and $t=3000$ only for the koopmons, thereby increasing the population on the ground state.

\end{document}